\documentclass[10pt, amsmath, amssymb]{amsart}

\usepackage{amsmath}
\usepackage{amsthm}
\usepackage{amssymb}
\usepackage{graphics}
\usepackage{bm}
\usepackage{natbib}

\newcommand{\be}{\begin{equation}}
\newcommand{\ee}{\end{equation}}
\newcommand{\ba}{\begin{eqnarray}}
\newcommand{\ea}{\end{eqnarray}}
\newcommand{\bi}{\begin{itemize}}
\newcommand{\ei}{\end{itemize}}
\newcommand{\bn}{\begin{enumerate}}
\newcommand{\en}{\end{enumerate}}
\newcommand{\bp}{\begin{proof}}
\newcommand{\ep}{\end{proof}}

\newcommand{\wt}{\ensuremath{\widetilde}}

\newcommand{\mr}{\ensuremath{\mathrm}}

\newcommand{\mc}{\ensuremath{\mathcal}}
\newcommand{\mf}{\ensuremath{\mathfrak}}
\newcommand{\ov}{\ensuremath{\overline}}
\newcommand{\sm}{\ensuremath{\setminus}}

\newcommand{\intfty}{\ensuremath{\int _{-\infty} ^{\infty}} }

\newcommand{\Om}{\ensuremath{\Omega}}

\newcommand{\La}{\ensuremath{\Lambda }}
\newcommand{\la}{\ensuremath{\lambda }}
\newcommand{\om}{\ensuremath{\omega}}
\newcommand{\ka}{\ensuremath{\kappa }}

\newcommand{\eps}{\ensuremath{\epsilon }}
\renewcommand{\bm}{\ensuremath{\mathbb }}

\newcommand{\ip}[2]{\ensuremath{\langle {#1} , {#2} \rangle}}

\newcommand{\dom}[1]{\ensuremath{\mathrm{Dom} ({#1}) }}
\renewcommand{\dim}[1]{\ensuremath{\mathrm{dim} \left( {#1} \right) }}

\newcommand{\ran}[1]{\ensuremath{\mathrm{Ran} ({#1}) }}
\renewcommand{\ker}[1]{\ensuremath{\mathrm{Ker} ({#1}) }}

\newcommand{\im}[1]{\ensuremath{\mathrm{Im} \left( {#1} \right) }}
\newcommand{\re}[1]{\ensuremath{\mathrm{Re} \left( {#1} \right) }}

\newcommand{\supp}[1]{\ensuremath{\mathrm{supp} ({#1}) }}

\numberwithin{equation}{section}

\newtheorem{thm}[subsubsection]{Theorem}
\newtheorem{claim}[subsubsection]{Claim}

\newtheorem{lemming}[subsubsection]{Lemma}

\newtheorem{cor}[subsubsection]{Corollary}

\begin{document}

\bibliographystyle{unsrt}

\title[Representation of symmetric operators in de Branges space]{Representation of simple symmetric operators with deficiency indices $(1,1)$ in de Branges space}

\author{R.T.W. Martin}

\address{Department of Mathematics \\ University of California, Berkeley\\
Berkeley, CA, 94720 \\
phone: +1 510 642 6550 \\ fax: +1 510 642 8204}

\email{rtwmartin@gmail.com}

\begin{abstract}

    Recently it has been shown that any regular simple symmetric operator with deficiency indices
$(1,1)$ is unitarily equivalent to the operator of multiplication in a reproducing
kernel Hilbert space of functions on the real line with the Kramer sampling property. This work
has been motivated, in part, by potential applications to signal processing and mathematical physics.
In this paper we exploit well-known results about de Branges-Rovnyak spaces and
characteristic functions of symmetric operators to prove that any such a symmetric
operator is in fact unitarily equivalent to multiplication by the independent variable
in a de Branges space of entire functions. This leads to simple new results on the spectra
of such symmetric operators, on when multiplication by $z$ is densely defined in de Branges-Rovnyak
spaces in the upper half plane, and to sufficient conditions for there
to be an isometry from a given subspace of $L^2 (\bm{R}, d\nu)$ onto a de Branges space of entire functions
which acts as multiplication by a measurable function.

\vspace{5mm}   \noindent {\it Key words and phrases}: symmetric (isometric) operators with deficiency indices $(1,1)$,
self-adjoint (unitary) extensions, de Branges spaces, de Branges-Rovnyak spaces, Lifschitz characteristic function,
Carath\'{e}odory or angular derivative.

\vspace{3mm}

\end{abstract}

\maketitle

\section{Introduction and Motivation}

    A reproducing kernel Hilbert space $\mc{H}$ of functions on a set
$X \subset \bm{C}$ is said to have the Kramer sampling property if it has a total
orthogonal set of point evaluation vectors \cite{Kramer}. The reason for this terminology
is clear: If $\delta _x$ denotes the point evaluation vector in $\mc{H}$ at
the point $x\in X$, \emph{i.e.} $\ip{f}{\delta _x} = f(x)$ for all $f \in \mc{H}$,
and if $\{ \delta _{x_n} \}$ is a total orthogonal set, then any $f \in \mc{H}$
can be reconstructed from its `samples', or values taken on the set of points $\{ x_n \} \subset X$,
using the sampling formula \be f = \sum _{n  } f(x_n)  \frac{\delta_{x_n}}{\delta _{x_n} (x_n)}.
\label{eq:recon} \ee We will only consider the case where $\mc{H}$ is separable, so that
all total orthogonal sets in $\mc{H}$ are countable.

    The classic examples of such spaces are the Paley-Wiener spaces $B(\Om )$ of $\Om$-bandlimited
functions, $\Om >0$. The space $B(\Om )$ is the image of $L^2 [-\Om , \Om]$ under the Fourier
transform, and for $f \in B(\Om)$ the reconstruction formula (\ref{eq:recon}) takes the form
\be f(x) = \sum _{n \in \bm{Z}}  f(x_n) \frac{ \sin \left( \Om (x-x_n ) \right) }{\Om (x-x_n)},
\label{eq:shan} \ee where $(x_n ) \subset \bm{R}$ is any sequence of points satisfying $x_{n+1} - x_n = \frac{\pi}{\Om}$. In this
particular case, the formula (\ref{eq:shan}) is called the Shannon sampling formula. These
spaces are used ubiquitously in signal processing to discretize and later reconstruct continuous signals;
\emph{e.g} audio signals to be recorded on compact disc, images for digital transmission, etc. The kernel of the idea is that by choosing
a sufficiently large bandlimit, a given continuous signal can be well-approximated by a bandlimited function.
The samples of this bandlimited approximation on a sufficiently dense set of points can then be recorded and
later used in Shannon's sampling formula to reconstruct the approximating function.

    There has been a significant amount of recent work focusing on the search for
and study of other reproducing kernel Hilbert spaces of functions on the real line with the sampling property
(or with the more general property that they contain a Riesz basis or frame of point evaluation vectors) in an
effort to develop more efficient methods of discretizing and reconstructing certain classes of
continuous signals in signal processing. See, for example \cite{Aldroubi} \cite{Everitt2} \cite{Silva} \cite{Kempfsamp}
\cite{Martin-symsamp} \cite{Garcia}, to name a few. Recently, it has been discovered that there is a connection between the
study of such spaces and the spectral representations of symmetric operators with deficiency
indices $(1,1)$. In particular, it has been shown that any regular simple symmetric operator, $B$,
with deficiency indices $(1,1)$ is unitarily equivalent to multiplication by the independent variable in a reproducing kernel
Hilbert space $\mc{H}$ with the sampling property where $\mc{H} \subset L^2 (\bm{R} , d\sigma)$ has the following properties
\cite{Silva} \cite{Kempfsamp} \cite{Martin-symsamp} \cite{Krein}.
First, the positive measure $\sigma$ can be chosen to be equivalent to Lebesgue measure, in fact $\sigma$ can
be chosen such that  $d\sigma (x) = \sigma ' (x)  dx $, where
$\sigma '  >0$ is continuous on $\bm{R}$. The space $\mc{H}$ consists of certain functions which
are meromorphic in $\bm{C}$ with no poles on the real axis. Furthermore, given any self-adjoint extension
$B ' $ of $B$, its spectrum $\sigma (B')$ is purely discrete and can be arranged as a monotonically strictly increasing
sequence of eigenvalues $(\la _n ) _{n \in \bm{Z}}$ with no finite
accumulation point, and if $\delta _\la $ denotes the point evaluation vector in $\mc{H}$ at any $\la \in \bm{R}$,
then $\{ \delta _{\la _n } \}$ is a total orthogonal set in $\mc{H}$.

    The main goal of this paper is to fully connect the theory developed in
\cite{Silva}  \cite{Kempfsamp} and \cite{Martin-symsamp} with the theory of de Branges spaces by showing that
any linear transformation $B$ with domain and range contained in a separable Hilbert space $
\mc{H}$ is regular, closed and simple symmetric with deficiency indices $(1,1)$ if and only if it is unitarily
equivalent to multiplication by the independent variable in a de Branges space of entire functions.
\footnote{ As observed in \cite{Silva} \cite{Kaltenback} if $B$ is actually an entire operator in the
sense of M.G. Krein \cite{Silva} \cite{Krein}, then it is already known that the space $\mc{H}$ described
in the previous paragraph is a de Branges space.} This will be accomplished by straightforwardly combining
known results about de Branges-Rovnyak spaces and Lifschitz' theory of the characteristic functions of
simple symmetric operators with deficiency indices $(1,1)$ \cite[Appendix 1]{Glazman}, \cite{Lifschitz}.

    In the process of achieving this result, it will be demonstrated that multiplication by
the independent variable is a simple symmetric linear transformation, with deficiency indices $(1,1)$
in any de Branges-Rovnyak space in the upper half plane provided this space is defined by an extreme point of the unit ball
in $H ^\infty$. New results on when this linear transformation is densely defined, and on the
spectra of simple symmetric operators with deficiency indices $(1,1)$ will also be presented. Finally,
new sufficient conditions for a subspace of $L^2 (\bm{R} , d\nu)$ to be the image of a de
Branges space of entire functions under an isometry which acts as multiplication by a measurable function
will be proven (see Theorem \ref{thm:RKHS} and Theorem \ref{thm:RKHS2} of Subsection \ref{subsection:subsamp}).

\subsection{Some notation}

     We will define \be \mu (z) := \frac{z-i}{z+i} ; \ \ \mu : \bm{U} \rightarrow \bm{D} \label{eq:CaT} \ee so that $\mu ^{-1} (z) = i
\frac{1+z}{1-z}$. Here $\bm{U}$ denotes the open complex upper half plane, $\bm{D}$ the unit disc, and $\bm{T}$, $\bm{L}$
will be used for the unit circle and open complex lower half plane, respectively.
Throughout this paper $\phi$ will denote a function in $B_1 \left( H^\infty ( \bm{U} ) \right)$, the
unit ball of $H^\infty (\bm{U})$, and
$\varphi := \phi \circ  \mu ^{-1}$ will be the canonical image of $\phi$ in $B_1 \left( H^\infty (\bm{D}) \right)$,
$\phi = \varphi  \circ \mu$. If $\phi$ is an extreme point of the unit ball of $H^\infty (\bm{U})$ we will
sometimes simply say that $\phi$ is extreme. Let $M$ denote the operator of multiplication by the independent variable in $L^2 (\bm{R})$,
so that $\phi (M)$ is multiplication by $\phi$. Let
$\mc{U} : H^2 (\bm{D}) \rightarrow H^2 (\bm{U} )$ denote the canonical isometry of $H^2 (\bm{D})$ onto $H^2 (\bm{U} )$.
If $f \in H^2 (\bm{D})$, \be \mc{U}  f (z) = \frac{1 -\mu (z)}{\sqrt{\pi} } f(\mu (z) ) \ \ \in H^2 (\bm{U}), \ee
and if $f \in H^2 (\bm{U})$, then  \be  \mc{U}  ^{-1} f(z) = \sqrt{\pi} \frac{f (\mu ^{-1} (z) )}{1-z}  \ \ \in H^2 (\bm{D} ).
\label{eq:isocanon} \ee  It follows that $\mc{U} ^{-1} \phi (M) \mc{U}$ acts as multiplication by $\varphi$ on $H^2 (\bm{D})$.

For $f \in L^\infty (\bm{R})$, let $T _f := P _{H^2 (\bm{U})} f(M) | _{H ^2 (\bm{U})} $, where $f(M)$ is
multiplication by $f$. The notation $K^2  _\phi$ will be used for the de Branges-Rovnyak
space which is the range of $R _\phi := \sqrt {1 - T_\phi T _{\ov{\phi}} }$, endowed with the inner product that makes $R _\phi$ a co-isometry of
$H^2 (\bm{U})$ onto its range.  Hence if $f,g \in H^2 (\bm{U})$ and at least one of them is orthogonal to $\ker{R _\phi}$, the
kernel of $R_\phi$, then $\ip{R_\phi f }{R_\phi g} _\phi = \ip{f}{g}$, where $\ip{\cdot}{\cdot} _\phi$ denotes the inner product in $K^2 _\phi$.
The notation $\mc{H} (E)$ is reserved for the de Branges space of entire functions defined
using the de Branges function $E$.

  It is not hard to check that $\mc{U}$ is an isometry
of $K^2 _\phi $ onto $K^2 _\varphi $, where $K^2 _\varphi $ is the usual de Branges-Rovnyak
space of functions in $H^2 (\bm{D})$, defined using $\varphi$, \emph{i.e.}, $K ^2 _
\varphi$ is the range of $\hat{R} _{\varphi}
:= \sqrt{1 - \hat{T} _{\varphi} \hat{T} _{\ov{\varphi}}}$ in $H^2 (\bm{D})$ endowed with the
inner product that makes $\hat{R} _{\varphi}$ a co-isometry of $H^2 (\bm{D} )$ onto $K^2 _{\varphi} $. Here $\hat{T} _\varphi$
is the Toeplitz operator $P_{H^2 (\bm{D})} \hat{M} _\varphi | _{H^2 (\bm{D})}$, and $\hat{M} _\varphi = \mc{U} ^{-1} \phi (M) \mc{U}$
acts as multiplication by $\varphi$ on $H^2 (\bm{D})$. The notation, $\ip{\cdot}{\cdot} _\phi$ and $\ip{\cdot}{\cdot} _\varphi$ will
be used for the inner products in $K^2 _\phi$ and $K^2 _\varphi$ respectively. Of course these inner products
reduce to the usual $L^2$ inner products in the case where $\phi = F$ is an inner function. If it is clear from the context
which space is being dealt with, the generic $\ip{\cdot}{\cdot}$ will sometimes be used. All inner products are assumed
to be conjugate linear in the second argument.

  For future reference, for $w \in \bm{U}$, we will let $k_w ^\phi$ denote the point evaluation vector at
$w$ in $K^2 _\phi$. That is $k_w ^\phi$ is the element of $K^2 _\phi$ such that $\ip{f}{k_w ^\phi} _\phi
= f(w)$ for all $f \in K^2 _\phi$. This vector has the form \be k_w ^\phi (z) := \frac{i}{2\pi}
\frac{1-\ov{\phi (w)} \phi (z)}{z- \ov{w}}. \ee Similarly, for $w \in \bm{D}$, $k_w ^\varphi$ will denote
the point evaluation vector in $K^2 _\varphi$ at $w$, given by the formula
\be k _w ^\varphi (z) := \frac{1-\ov{\varphi (w)}\varphi (z)}{1- \ov{w} z} .\ee Finally for any function
$f$, analytic on a region $\Om$, $f^*$ will denote the function $f^*(z) = \ov{f (\ov{z} )}$ analytic
in the reflection of the region $\Om$ in the real axis.

\section{Representation of an arbitrary simple symmetric linear transformation with deficiency indices $(1,1)$}

    Given a Hilbert space $\mc{H}$, let $V$ denote an arbitrary closed simple isometric linear transformation
with deficiency indices $(1,1)$ in $\mc{H}$. Here, recall that the deficiency indices $(n_+ , n_-)$ of $V$ are
defined as $n_+ := \dim{ \dom{V}  ^\perp} $ and $n_- := \dim{\ran{V} ^\perp}$, and that an isometric
linear transformation is called simple if it has no unitary restriction to any non-trivial subspace. Further recall that
a linear transformation in $\mc{H}$ is called closed if its graph is a closed subspace of $\mc{H} \oplus
\mc{H}$. An isometric linear transformation $V$ is bounded, and so will be closed provided its domain is a
closed subspace of $\mc{H}$.

    Given such a $V$, let $B : = \mu ^{-1} (V) = i (1+V) (1-V) ^{-1}$, where $(1-V) ^{-1}$ is defined as a
 linear map from $\ran{V-1} $ onto $\dom{V}$. The map $B$ is a closed simple symmetric linear
transformation in $\mc{H}$ with deficiency indices $(1,1)$. Recall that a symmetric linear transformation is called
simple if it has no self-adjoint restriction to a dense domain in a non-trivial subspace, and that the
deficiency indices of $B$ are also $(n_+ , n_- )$ where $n _+ = \dim{\dom{V} ^\perp} = \dim{\ran{B+i} ^\perp} $
and $n_- = \dim{\ran{V} ^\perp} = \dim{\ran{B-i} ^\perp}$. Note that if $B$ is densely defined, then its
adjoint $B^*$ is a well-defined closed operator and $\ran{B-\ov{z}} ^\perp = \ker{B^* -z}$ for all $z \in \bm{C}$.
Finally recall that $V = \mu (B)$, and that the map
$B \mapsto \mu (B)$ is a bijection of the set of simple symmetric linear transformations with deficiency
indices $(1,1)$ onto the set of simple isometric linear transformations with deficiency indices $(1,1)$. If $B$
is symmetric, $V := \mu (B)$ is called its Cayley transform. For the basic theory on symmetric operators, their
Cayley transforms and deficiency indices, see for example \cite{Glazman} or \cite{Reed2}.

    The Lifschitz characteristic function $w_V$ of a simple isometric linear transformation $V$ with
deficiency indices $(1,1)$ is defined as
\be w _V (z) := \frac {z \ip{ (U-z) ^{-1} \psi _+}{\psi _+}}{\ip{ (U-z) ^{-1} U \psi _+}{\psi _+} }, \label{eq:carv} \ee
where $U$ is an arbitrary unitary extension of $V$, and $ 0
\neq \psi _+ \in \dom{V} ^\perp$ \cite[Appendix 1]{Glazman} \cite{Lifschitz}. Note that if
$B$ is defined as $B:= \mu ^{-1} (V)$, then $\dom{V} ^\perp = \ran{B +i} ^\perp$. Choosing a different
unitary extension $U'$ in the definition of $w _V$ only changes $w _V$ by a unimodular constant. For
this reason $w_V$ is defined only up to such a constant, and we will say any two characteristic
functions coincide if they differ by such a constant. The function $w_V$ belongs to $B_1 \left( H^\infty
(\bm{D} ) \right)$, and obeys $w_V (0) = 0$.

The Lifschitz characteristic function of
a simple symmetric linear transformation with deficiency indices $(1,1)$, $B$, is then defined
as $\om _B ( \la ) := w _V (\mu (\la) )$. A short calculation shows that for $a, \la \in \bm{C} \sm \{-i \}$,
\ba \left( \mu (a) -\mu (\la ) \right) ^{-1} & = & \left[ \frac{1}{\la +i} \ \frac{a-i}{a+i}
\left( (\la +i) -(\la -i) \frac{a+i}{a-i} \right) \right] ^{-1} \nonumber \\
& = & (\la +i) \frac{a+i}{a-i} \left(  2i \frac{a-\la}{a-i} \right) ^{-1} \nonumber \\
& = & \frac{\la +i}{2i} \frac{a+i}{a-\la} \nonumber \\ & = & \frac{\la +i}{2i} \frac{a+i}{a-i} \left( 1 + \frac{\la -i}{a-\la}
\right). \ea If $U$ is a unitary extension of $V =\mu (B)$, then  $A:= \mu ^{-1} (U)$ is a
self-adjoint extension of $B$ and
\be (U-\mu (\la) )^{-1} = (\mu (A) - \mu (\la ) ) ^{-1} =
\frac{\la +i}{2i} \mu(A) ^* \left( 1 + (\la -i) (A-\la ) ^{-1}
\right). \ee It follows that the characteristic
function of $B$ can be written as
\be \om _B (\la ) = w _{\mu (B)} ( \mu (\la ) ) = \frac{\la -i}{\la +i} \frac{\ip{ ( I + (\la -i)
(A-\la ) ^{-1} ) \psi _+ }{ \mu(A) \psi _+ } }{\ip{(I + (\la -i) (A-\la) ^{-1} ) \psi _+ }{\psi _+} }. \label{eq:symchar} \ee
Again, $\om _B$ is defined modulo unimodular constants, $\om _B \in B_1 \left( H^\infty (\bm{U} ) \right)$, and $\om _B (i) =0$.
Alternatively, one can calculate that $\left( \mu (a) -\mu (\la ) \right) ^{-1} = \dfrac{\lambda
 +i}{2i}\left(1+(\lambda+i)(a-\lambda)^{-1}\right)$, so that the characteristic function can also be
 written as \ba \omega_B(\lambda) & =&
 \frac{\lambda-i}{\lambda+i}
 \frac{\langle \left(I+(\lambda+i)(A-\lambda)^{-1}\right)
 \psi_+,\psi_+\rangle}{\langle
 \left(I+(\lambda-i)(A-\lambda)^{-1}\right)
\psi_+,\psi_+\rangle}.
 \ea

\subsubsection{Remark} \label{subsubsection:technik} There is a slight technicality, glossed over in the above discussion that should be pointed out. It is straightforward
to show that $B$ is densely defined if and only if no unitary extension $U$ of the Cayley transform $V = \mu (B)$ of $B$ has $1$ as an eigenvalue
(see the last part of the proof of Theorem \ref{thm:noone} to come). Hence if $B$ is not densely defined, then there is a unitary extension $U _1$
of $V$ which has $1$ as an eigenvalue so that for this particular $U_1$, $A = \mu ^{-1} (U_1)$ is not well defined. This unitary extension
is unique, if any other unitary extension $U'$ also has $1$ as an eigenvalue then one can calculate that $1$ would have to be an eigenvalue
of $V$, contradicting the simplicity of $V$. In conclusion, if $B$ is not densely defined, then in the formula (\ref{eq:symchar}) for the characteristic
function of $B$, $A = \mu ^{-1} (U)$ must be chosen as the inverse Cayley transform of a unitary extension $U$ of $V=\mu (B)$ where $U\neq U_1$.

    We will employ the following theorems of M. S. Lifschitz on such isometric and symmetric linear transformations
and their characteristic functions \cite[Appendix 1]{Glazman} \cite{Lifschitz}:

\begin{thm}{ (Lifschitz) }
In order that two simple isometric linear transformations with deficiency
indices $(1,1)$ be unitarily equivalent, it is necessary and sufficient that their characteristic
functions coincide. \label{thm:charequiv}
\end{thm}

\begin{thm}{ (Lifschitz) }
    Given any $w \in B_1 (H^\infty (\bm{D}) )$  obeying $w(0) =0$,  there
exists a simple isometric linear transformation with deficiency indices $(1,1)$
whose characteristic function is $w$. \label{thm:exist}
\end{thm}

\subsubsection{Remark} \label{rmk:bi} Using the bijective correspondence between simple isometric and symmetric
linear transformations with deficiency indices $(1,1)$, the analogous statements obtained by
replacing $w \in B_1 \left( H^\infty (\bm{D} ) \right)$ with $\om \in B_1 \left( H^\infty (\bm{U} ) \right)$,
isometric with symmetric, and $w(0) = 0$ with $\om (i) = 0 $ in the above two theorems are also true. \vspace{.05truein}

    Let $\varphi \in B_1 (H^\infty (\bm{D} ))$, and consider the function
$b _\varphi:= \frac{1+\varphi}{1-\varphi}$. Then $b_\varphi$ has non-negative real part, and so has the Herglotz integral
representation
\be b_\varphi (z) = \int _0 ^{2\pi} \frac{e^{i\theta} +z}{e^{i\theta} -z} d \rho _\varphi (e ^{i\theta}) + i \im{b_\varphi (0)},
\label{eq:Herglotz} \ee
where $\rho _\varphi $ is a positive, finite Borel measure on the unit circle $\bm{T}$ and $\im{z}$ denotes the imaginary
part of z.
Now consider the space $L^2 _\varphi (\bm{T} )$ of functions on the unit circle square integrable with respect
to this measure $\rho _\varphi$, and let $Z _\varphi$ denote the unitary operator of multiplication by $z$ in $L^2 _\varphi (\bm{T})$,
\emph{i.e.} $Z _\varphi f(e^{i\theta} ) = e^{i\theta} f (e^{i\theta} )$. Let $\psi _+ \in L^2 _\varphi (\bm{T} )$ be defined by $\psi _+ (z) := 1/z$.
If $\dom{Z ' _\varphi} := \{ f \in L^2 _\varphi (\bm{T} ) | \  \ip{f}{\phi _+ } =0  \}$, then it is not hard to check that
$Z '_\varphi := Z _\varphi |_{\dom{Z ' _\varphi}}$ is a closed simple isometric linear transformation with deficiency indices $(1,1)$ whose
characteristic function is \be w _\varphi (z) = \frac{\varphi(z) - \varphi(0)}{1 -\ov{\varphi(0)} \varphi (z)}
\label{eq:charfunction} .\ee
In particular, if $\varphi (0) = 0$ then $w _\varphi (z) = \varphi (z)$. These facts are collected in the following lemma.

\begin{lemming}
    Let $V$ be a simple isometric linear transformation with deficiency indices $(1,1)$ and
characteristic function $w$. Then $V$ is unitarily equivalent to $Z ' _w: \dom{Z ' _w} \rightarrow L^2 _w (\bm{T} )$,
where $Z '_w$ acts as multiplication by $z$ on its domain.
If $\varphi \in B_1 (H^\infty (\bm{D} ) )$, then the characteristic
function $w_\varphi$ of $Z' _\varphi$ is given by equation (\ref{eq:charfunction}). \label{lemming:charfunction}
\end{lemming}

    The following proof of this lemma follows immediately from that of Theorem \ref{thm:exist}
in \cite[Appendix 1]{Glazman}.

\begin{proof}
    The final statement will be proven first. Suppose that $\varphi \in B_1 (H^\infty (\bm{D})$.
By definition (see (\ref{eq:carv})), the characteristic function $v(z)$ of $Z' _\varphi$ is equal to
\be  v(z) = \frac{z \ip{(Z_\varphi -z) ^{-1} \psi _+ }{\psi _+}}{\ip{ (Z_\varphi -z) ^{-1} Z_\varphi \psi _+}{\psi _+}}, \ee
where $\psi _+ (z) = 1/z$ in $L^2 _\varphi (\bm{T})$. First observe that
\be \ip {(Z_\varphi - z) ^{-1} \psi _+ }{\psi _+ } = \int _0 ^{2\pi} \frac{1}{e^{i\theta} -z} d\rho _\varphi (e^{i\theta}),\ee
while \be   \ip{(Z_\varphi - z) ^{-1} Z_\varphi \psi _+ }{\psi _+ } = \int _0 ^{2\pi} \frac{e^{i\theta} }{e^{i\theta} -z} d\rho _\varphi (e^{i\theta}). \ee
By the Herglotz representation (\ref{eq:Herglotz}) of $b_\varphi$, $\rho _\varphi ( \{ \bm{T} \} ) = \re{b_\varphi(0)}$ and
\ba b_\varphi (z) - b_\varphi (0) & = & \int _0 ^{2\pi} \frac{e^{i\theta} +z }{e^{i\theta} -z} d\rho _\varphi (e^{i\theta}) - \re{b_\varphi (0)} \nonumber \\
&= & \int _0 ^{2\pi} \frac{e^{i\theta} +z - (e^{i\theta} -z ) }{e^{i\theta} -z} d\rho _\varphi (e^{i\theta}) \nonumber \\
& = & 2z \ip{(Z _\varphi -z) ^{-1} \psi _+ }{\psi _+}. \ea Similarly,
\be b_\varphi (z) + \ov{b_\varphi (0)} = 2 \ip{(Z _\varphi -z) ^{-1} Z_\varphi \psi _+ }{\psi _+}.  \ee Using that
$b_\varphi = \frac{1+\varphi}{1-\varphi}$, we conclude that
\be v = \frac{b_\varphi  - b(0)}{b_\varphi + \ov{b(0)}} = \chi \frac{\varphi - \varphi(0)}{1 - \ov{\varphi(0)} \varphi }, \ee
where $\chi = \frac{1- \ov{\varphi(0)}}{1-\varphi(0)} \in \bm{T}$. This establishes the formula (\ref{eq:charfunction}),
and the final statement of the lemma.

    To prove the remainder of the lemma, recall that if $w$ is the characteristic function of $V$, then
$w(0) = 0$. By (\ref{eq:charfunction}), the characteristic function of $Z_w '$ is also $w$. The rest of the lemma now follows
from Theorem \ref{thm:charequiv}.
\end{proof}

\subsection{An isometry of $K^2 _\varphi$ onto $H^2 _\varphi$}

    Let $\rho _\varphi$ be the positive Borel measure on the circle
$\bm{T}$ defined by the function $\varphi \in B_1 (H^\infty
(\bm{D} ) )$ as in equation (\ref{eq:Herglotz}). The Cauchy integral of $\rho _\varphi$ is defined by
\be \Gamma _{\varphi} (z) = \int _\bm{T} \frac{1}{1-\ov{w} z} d \rho _\varphi (w),\ee
\cite[Chapter III]{Sarason-dB}. Given $f \in L^2 _\varphi (\bm{T})$, one defines
$\Gamma _\varphi f := \int _\bm{T} \frac{f(w)}{1-\ov{w}z} d\rho _\varphi (w) $. Let $H^2 _\varphi $ denote the closure
of the polynomials in $L^2 _\varphi$. For any $f \in H^2 _\varphi $, define
$W_\varphi f (z) := (1 - \varphi (z) ) \Gamma _\varphi f (z)$. Then
as shown, for example in \cite[III-7]{Sarason-dB}, $W_\varphi$ is
an isometry of $H^2 _\varphi $ onto $K^2 _\varphi$.

    It follows from \cite[III-8]{Sarason-dB} that \be W_\varphi Z _\varphi ^* W_\varphi ^{-1}
= X + (1 - \varphi (0) ) ^{-1} \ip{\cdot}{k_0 ^\varphi} _\varphi S^*\varphi, \label{eq:shifty} \ee
where $S^*$ is the backward shift in $H^2 (\bm{D} )$ and $X := S^* | _{K^2 _\varphi}$. Furthermore,
as shown in \cite[III-7]{Sarason-dB}, the image of the constant function $f(z) =1$
under $W _\varphi$ is $\left( 1- \varphi (0) \right) ^{-1} k _0 ^\varphi$.

\subsection{The case where $\varphi $ is extreme}

    Recall that $\varphi $ is an extreme point of the unit ball if and only if $\ln (1 - | \varphi | )$
is not integrable \cite[pg. 138]{Hoffman}. If $H^2 _\varphi $ is the closure of the
polynomials in $L^2 _\varphi $, further recall that $H^2 _\varphi = L^2 _\varphi $ if and only if the Radon-Nikodym
derivative of the absolutely continuous part of $\rho _\varphi $ with respect to Lebesgue measure is not
log-integrable \cite[pg. 50]{Hoffman}. This derivative is equal to $\frac{1 -| \varphi | ^2}{|1 - \varphi | ^2}$
and hence is not log-integrable if and only if $\varphi $ is an extreme point. \vspace{.05truein} \\

   Since $Z _\varphi ^*$ agrees with $(Z' _\varphi) ^*$ on $\ran{Z' _\varphi}$, $\ran{Z' _\varphi}$ is the orthogonal complement of
the constant function $1$ in $L^2 _\varphi$, and the image of $1$ under the isometry $W_\varphi$
is a constant multiple of $k _0 ^\varphi$, it follows that the image of $\ran{Z' _\varphi}$ under
the isometry $W_\varphi$ is the subspace $\mc{S} _0 \subset K^2 _\varphi$ of functions in $K^2 _\varphi$ which vanish at $0$. It then
follows from equation (\ref{eq:shifty}) that the image of $(Z' _\varphi) ^*$ under $W_\varphi$ is $S^* | _{\mc{S} _0}$.

If $f(z) = z g(z) \in \mc{S} _0$, then
it is clear that $g = S^* f \in K^2 _\varphi$ since $K^2 _\varphi$ is invariant under $S^*$. Moreover,
by \cite[II-9]{Sarason-dB}, with $X = S^* | _{K^2 _\varphi}$,
\be X^*Xf = X^* S^* f = SS^*f - \ip{S^*f}{S^*\varphi} _\varphi \varphi = f - \ip{S^*f}{S^* \varphi} _\varphi \varphi.
\label{eq:shiftiso} \ee
Since $\varphi$ is an extreme point, it follows from \cite[V-3]{Sarason-dB},
that $\varphi \notin K^2 _\varphi$. Hence it
must be that the inner product $\ip{S^*f}{S^* \varphi} _\varphi =0$, that $X^* X f = S S^* f = f$ and that $X^* | _{S^* \mc{S} _0}
= S | _{S^* \mc{S} _0} $. It follows that $S | _{S^* \mc{S} _0} =W_\varphi Z'_\varphi W_\varphi ^{-1}$ is an isometry from
the subspace $S^* \mc{S} _0 $ onto $ \mc{S} _0$.

   The image of $S$ under the canonical isometry, $\mc{U}$, of $H^2 (\bm{D} )$ onto $H^2 (\bm{U})$ is $\mu (M)$, which
acts as multiplication by $\mu (z) = \frac{z-i}{z+i}$. The image of $\mc{S} _0$ under $\mc{U}$ is the subspace
$\mc{S} _i \subset K^2 _\phi$ of functions in $K^2 _\phi$ which vanish at $z=i$. It follows that if
$V ^\phi := \mc{U} S^* \mc{U} ^{-1} | _{\mc{S} _i} = \mc{U} W_\varphi Z' _\varphi W_\varphi ^* \mc{U}^* | _{\mc{S} _i = \mc{U} W_\varphi
\dom{Z ' _\varphi} }$, then $(V ^\phi) ^* = 1/ \mu (M) | _{\mc{S} _i}$ where $1/\mu(M)$ is multiplication by $1/\mu$,
and $V ^\phi = \mu(M)  | _{\ran{(V _\phi) ^*}}$ is the image of $Z' _\varphi$ under the isometry $\mc{U}W_\varphi$.
Define the linear transformation $(V ^\phi -1)^{-1}: \ran{V ^\phi -1} \rightarrow \dom{V ^\phi}$. Then the inverse Cayley transform of $V ^\phi$,
$M ^\phi := \mu ^{-1} (V ^\phi) = i (1+V ^\phi) (1-V ^\phi) ^{-1}$, is a symmetric linear transformation which acts as multiplication by $z$ on the domain
$\dom{M ^\phi} := \ran{V ^\phi -1} \subset K^2 _\phi$, $M^\phi = M | _{\dom{M^\phi}}$.
This domain will be dense if and only if $\ran{V ^\phi -1}$ is dense in $K^2 _\phi$.
We will provide necessary and sufficient conditions on $\phi$ for $M ^\phi$ to be a densely defined symmetric operator
in $K^2 _\phi$ in the case where $\phi$ is extreme in Section \ref{section:dense}.

    Applying Lemma \ref{lemming:charfunction} and Remark \ref{rmk:bi} immediately yields:

\begin{thm}
    The transformation $V ^\phi := \mu(M) | _{1/ \mu(M) \mc{S} _i }$ is a closed simple isometric transformation
with deficiency indices $(1,1)$ in $K^2 _\phi$. Its inverse Cayley transform, $M ^\phi = \mu ^{-1} (V ^\phi )$
is a closed simple symmetric linear transformation with
deficiency indices $(1,1)$ which acts as multiplication by $z$ on the domain $\dom{M ^\phi} =
\ran{V ^\phi -1} \subset K^2 _\phi$. The characteristic function of $M ^\phi$ is  \be \om _\phi =
\frac{\phi - \phi (i)}{ 1 - \ov{\phi (i)} \phi } = w _\varphi \circ \mu . \label{eq:multchar} \ee

    Conversely, if $B$ is any closed simple symmetric linear transformation
with deficiency indices $(1,1)$ whose characteristic function $\om $ is an extreme
point, then $B$ is unitarily equivalent to $M ^{\phi _\alpha}$ in the de Branges-Rovnyak space
$K^2 _{\phi _\alpha}$ where $\phi _\alpha  = \frac{ \om - \alpha}{ 1 - \ov{\alpha} \om } $ and
$\alpha = -\phi _\alpha (i) \in \bm{D}$ is arbitrary. \label{thm:symrep}
\end{thm}

    Again, as observed in Remark \ref{rmk:bi}, an analogous statement to the second part
of the above theorem holds for any simple isometric linear transformation with indices
$(1,1)$.

\subsubsection{Remark} In particular, choosing $\alpha =0$ in the above theorem shows that $B$ is
unitarily equivalent to multiplication by $z$ in $K^2 _\om $. It is easy to check that $\phi $ is inner or extreme
if and only if $\om $ is.

\subsubsection{Remark} If $M$ denotes multiplication by $x$ in $L^2 (\bm{R})$ and $F$ is an inner
function, it is not difficult to verify that $M^F$ is the unique closed symmetric restriction of $M$
to a linear subspace of $K^2 _F$ with deficiency indices $(1,1)$. \label{subsubsection:unique}

\subsubsection{A conjugation which commutes with $M ^\phi$} In the case where
$\varphi$ is extreme, so that $W _\varphi$ is an isometry of $L^2 _\varphi $
onto $K^2 _\varphi$, let $C$ denote the conjugation defined on
$L^2 _\varphi $ by $Cf(z) = 1/z \ov{f(z)}$ (see for example \cite{Lotto}). Recall that a conjugation is
an idempotent norm-preserving anti-linear map. It is easy to observe
that $CZ _\varphi = Z^* _\varphi C$, that $C \dom{Z ' _\varphi } = \ran{Z' _\varphi }$, and that $C\dom{Z' _\varphi}
^\perp = \ran{Z' _\varphi } ^\perp$. With these facts one can show that if $C _\phi$ denotes the image of $C$
in $K^2 _\phi$, $C _\phi := \mc{U}W_\varphi C W_\varphi ^* \mc{U} ^*$, that $C_\phi k _\la  ^\phi (z) = \frac{1}{2\pi i }
\frac{\phi (z) - \phi (i) } {z-i} $, and that the following lemma holds.
\label{subsubsection:conju}
\begin{lemming}
    The conjugation $C _\phi$ maps $\dom{M ^\phi}$ onto itself and commutes with $M ^\phi$. For
any $\la \in \bm{C}\sm\bm{R} $, $C_\phi \ran{M ^\phi -\la} ^\perp = \ran{M ^\phi - \ov{\la} } ^\perp$.
\label{lemming:con}
\end{lemming}

Any $\phi \in B_1 (H^\infty (\bm{U} ) )$ has the canonical representation
\be \phi(z) = \chi \prod _{n=1} ^\infty \frac{1 - z/ z_n}{1 - z / \ov{z_n}} \
e^{i\sigma z} \exp \left( i \intfty \frac{1+tz}{t-z} d\nu (t) \right) \label{eq:canon},\ee
where $(z_n ) _{n \in \bm{N}} \subset \bm{U}$, $\sum _n \left| \im{\frac{1}{ z_n }}  \right|
< \infty$, $\sigma \geq 0$, $\chi \in \bm{T}$ and $\nu$ is a finite, positive Borel measure on $\bm{R}$.
This formula actually defines a function which is analytic everywhere
in the region $\Om$ where $\Om := \bm{C} \sm \supp{\phi} ^*$ and $\supp{\phi} ^* $ is defined
as the union of the closure of the set $\{ \ov{z_n} \} _{n \in \bm{N}}$ with the closed
support of the measure $\nu$ on $\bm{R}$. Recalling that $\phi ^*$ is defined by
$\phi ^* (z) = \ov{\phi (\ov {z} ) } $, equation (\ref{eq:canon})
implies that $\phi(z) \phi ^* (z) = 1 $ for all $z \in \Om$.

    It is clear that if $h \in H^2 (\bm{U})$ then $h^* \in H^2 (\bm{L})$, and visa versa. If $\phi =F$
is inner, so that $F(x) \ov{F(x)} =1 $ almost everywhere $x \in \bm{R}$, then the non-tangential
limits of $F$ as $z$ approaches the real axis from either the upper or lower half-planes are
equal to $F(x)$ almost everywhere. To see this, note that if $z \in \bm{L}$, $F(z) = 1/ F^* (z)$,
and the non-tangential limits of $F (z)$ as $z$ approaches the real axis from below equal
$ 1 / \ov{F(x)} = F(x)$ almost everywhere. Hence if $h \in H^2 (\bm{U})$, $F h^*$ is that function whose non-tangential
limits as $z$ approaches the real axis from below equal $F(x) \ov{h (x)}$ almost everywhere, and
if $f,g \in H^2 (\bm{U})$ then $\ip{Ff}{g}=\ip{f}{F^*g}$.  Let $*$ denote the
operation $f \mapsto f^*$, and $F(M)$ be multiplication by $F$.

\begin{claim}
    If $F$ is inner, the mapping $\wt{C} _F := F(M) \circ *$ is a conjugation on $K^2 _F
= H^2 (\bm{U})\ominus F H^2 (\bm{U})$, and $f \in K^2 _F$ if and only if both $f, F f^*$ belong to $H^2(\bm{U})$.
\label{claim:inspace}
\end{claim}

\begin{proof}
  Suppose $f \in K^2 _F$. Then if $h \in H^2 (\bm{U})$ is arbitrary, $\ip{F^* f}{h}
=\ip{f}{Fh} =0$ so that $F^*f \in L^2 (\bm{R} ) \ominus H^2 (\bm{U})= H^2 (\bm{L})$ and $\wt{C}_F f = Ff^* \in H^2 (\bm{U})$.
Conversely, suppose that $f, Ff^* \in H^2(\bm{U})$. Then given any $Fh \in FH^2(\bm{U})$, $\ip{f}{Fh}
=\ip{F^*f}{h} = 0$ since $F^*f \in L^2 (\bm{R}) \ominus H^2(\bm{U})$. It follows that $f \in H^2 (\bm{U})\ominus
FH^2 (\bm{U})= K^2 _F$.

    The above shows that if $f \in K^2 _F$, then so is $\wt{C} _F f $. The mapping
$\wt{C} _F$ is clearly norm-preserving and anti-linear, so to show it is a conjugation,
it remains to show that it is idempotent. If $f \in K^2 _F$, then
$\wt{C} _F ^2 f = F (F f^*)^* = F F^* f = f$.
\end{proof}

    More generally if $\phi \in B_1 (H^\infty)$ is an extreme point, let
$\wt{C} _\phi = \phi (M) \circ *$.

\begin{claim}
    The map $\wt{C} _\phi$ defined on the linear span of the point evaluation vectors $k_w ^\phi$,
$w \in \bm{U}$ is anti-linear, norm-preserving and idempotent with range contained in $K^2 _\phi$. Its
unique continuous and anti-linear extension to all of $K^2_\phi$ is the conjugation $C _\phi$.
\end{claim}

\begin{proof}
    Consider the dense linear manifold in $K^2 _\phi$ which consists of all finite linear combinations
of the point evaluation vectors $k_w ^\phi $, for $w \in \bm{U}$.  It follows that
$\wt{C} _\phi k_w ^\phi (z) = \phi(z) (k_w ^\phi ) ^* (z) = \phi(z) \frac{i}{2\pi} \frac{1- \phi (w) \phi ^* (z)}{w-z} =
\frac{-i}{2\pi}\frac{\phi (z) - \phi (w)}{z-w}$. As mentioned in \ref{subsubsection:conju} above, this is equal to
$C_\phi k_w ^\phi$, so that $\wt{C} _\phi $ and $C_\phi$ are norm preserving, and anti-linear maps that
agree on a dense linear subspace of $K^2 _\phi$. It follows that $\wt{C} _\phi$ can be extended by anti-linearity
and continuity to all of $K^2 _\phi$ and that this extension agrees with $C _\phi$.
\end{proof}

\subsubsection{The characteristic function of $M ^\phi$}
When $\phi$ is extreme, one can also directly calculate the characteristic function
$\om _\phi = \om _{M ^{\phi}}$ of $M^\phi$ as follows. Consider the formula for the characteristic function
$\om _B$ of a simple symmetric linear transformation $B$ with deficiency indices $(1,1)$ as given
in equation (\ref{eq:symchar}),
\be \om _B (\la ) := w _V ( \mu (\la ) ) = \frac{\la -i}{\la +i} \frac{\ip{ ( I + (\la -i)
(A-\la ) ^{-1} ) \psi _+ }{ \mu (A) \psi _+ } }{\ip{(I + (\la -i) (A-\la) ^{-1} ) \psi _+ }{\psi _+} }. \ee

If $U$ is an arbitrary unitary extension of $\mu (B)$ which does not have $1$ as an eigenvalue (see Remark
\ref{subsubsection:technik}), and $A = \mu ^{-1} (U)$, it is not
difficult to show that the operator $ (A - z)(A - z') ^{-1} =
 \left( I + (z-z') (A -z )^{-1} \right)$ maps $\ran{B -\ov{z'}} ^\perp $ onto $\ran{B -\ov{z}} ^\perp$
(\cite{Krein}, pg. 9). Consider $B= M ^\phi$, the symmetric transformation of multiplication by $z$
in $K^2 _\phi$, where $\phi$ is an extreme point of the unit ball in $H^\infty$.
Observe that $k_z ^\phi \in \ran{M ^\phi - z } ^\perp$
for all $z \in \bm{U}$. Indeed, for any $f \in \dom{M ^\phi }$, $\ip{(M ^\phi -z) f}{k_z ^\phi} _\phi =0$. Note
that if $M ^\phi$ is actually densely defined, then $\ran{M^\phi -z} ^\perp = \ker{(M ^\phi )^* -\ov{z} }$ so
that each $k_z ^\phi$ is an eigenvector to the adjoint of $M^\phi$ with eigenvalue $\ov{z}$.
Let $\om  _\phi = \om _{M ^\phi}$, be defined using an arbitrary unitary extension $U$ of the Cayley transform of $M ^\phi$,
and let $M' = A$ be the inverse Cayley transform
of $U$. Since $\left( I + (z-z') (M' -z )^{-1} \right)
\psi _+ \in \ran{M - \ov{z} } ^\perp$, $\ran{M - \ov{\la}} ^\perp = \bm{C} \{ k_\la \} $, the one dimensional
subspace spanned by $k_\la$, and
$\ran{M -\ov{\la} } ^\perp = \bm{C} \{ C_\phi k_\la \}$ for each $\la \in \bm{U}$ by Lemma \ref{lemming:con},
it follows that $\left( I + (\la - i) (M' -\la )^{-1} \right)
\psi _+ = c(\la ) C_\phi k_\la$. Here $c(\la) \in \bm{C}$. Moreover, we can choose $\psi _- = k_i / \| k_i \|$ and
$\psi _+ = C _\phi k_i / \| k _i \|$. It follows that
\ba \om  _\phi (\la) & = & \frac{\la -i}{\la+i} \frac{\ip{c(\la) C_\phi k_\la}{k_i}}{\ip{c(\la) C_\phi k_\la}{C_\phi k_i}}
\nonumber \\ & =& \frac{\la -i}{\la +i} \frac{\ip{C_\phi k_\la}{k_i}}{\ip{k_i}{k_\la}}  =  \frac{\la -i}{\la +i} \frac{ (C_\phi k_\la ) (i) }{ k_i (\la)} \nonumber \\
& = & \frac{ \phi (i) - \phi (\la )}{ 1 - \ov{\phi (i)} \phi (\la )}, \label{eq:charfunk}  \ea
which is the same (up to a unimodular constant) as the formula for $\om _\phi = \om _{M ^\phi}$
given in the statement of Theorem \ref{thm:symrep}.

\section{Necessary and sufficient conditions for $M ^\phi$ to be densely defined.}
\label{section:dense}

    Given a simple symmetric linear transformation $B$ with deficiency indices $(1,1)$,
the following theorem of Lifschitz provides a necessary and sufficient condition on
the characteristic function $\om $ of $B$ for $B$ to be densely defined \cite[Appendix 1]{Glazman}:

\begin{thm}{ (Lifschitz) }
    $B$ will be densely defined if and only if $\lim _{\la \rightarrow \infty}
\la (\om (\la ) - e^{i\alpha} ) = \infty $; $0 < \eps \leq \arg{\la} \leq \pi - \eps$,
for each $\alpha \in [0 , 2\pi)$ and any fixed $\eps >0$. \label{thm:charfun}
\end{thm}

    Let $w:= \om \circ \mu ^{-1} $ be the characteristic function of $V := \mu (B)$. It
is not hard to see that the above necessary and sufficient condition for $B$ to be
a densely defined operator is related to the existence of the angular derivative of $w$
at the point $z=1$.

\subsection{Angular derivatives}

    Given $b \in B_1 \left( H^\infty (\bm{D} ) \right)$, $b$ has the canonical factorization
\be b(z) = \chi z^n \prod _n \frac{|a_n | }{a_n} \frac{a_n -z}{1-\ov{a_n} z} \label{eq:dcanon} \exp
\left( - \int _\bm{T} \frac{\zeta +z}{\zeta -z} d\rho (\zeta ) \right), \ee where
$\sum _n (1 - |a_n |) < \infty$, $\chi \in \bm{T}$ and $\rho $ is a positive Borel measure on $\bm{T}$.

    The function $b$ is said to have an angular derivative at $z \in \bm{T}$ in
the sense of Carath\'{e}odory if the non-tangential limits of $b, b'$ exist at $z$ and
the non-tangential limit of $b$ at $z$ has unit modulus. Define \be A _b := \left\{ z \in \bm{T} \ \left|  \ \sum _n \frac{1-|a_n|^2}{|z-a_n| ^2}
+ \int _\bm{T} |\theta -z| ^{-2} d \rho (\theta ) < \infty \right. \right\}. \ee The
following theorem, taken from \cite{Fricain}, is a combination of results of
\cite{Clark1} and \cite{Sarason-dB}.

\begin{thm}{ (Carath\'{e}odory, Ahern and Clark, Sarason) }
   Given $b \in B_1 (H^\infty (\bm{D} ) )$, and $\gamma \in \bm{T}$, the following
are equivalent: \\
(i) There exists a unimodular constant $c$ such that $\frac{b(z) - c}{z- \gamma}
\in K^2 _b$. \\
(ii) $\gamma \in A_b$. \\
(iii) $\lim \inf _{z \rightarrow \gamma} \frac{1 - |b(z)|}{1 -|z|} < \infty $. \\
(iv) $b$ has an angular derivative in the sense of Carath\'{e}odory at $\gamma$.\\
(v) every $f \in K^2 _b$ has a non-tangential limit at $\gamma$. \\
    Furthermore, if any of the above conditions hold, $c = b(\gamma) := \lim _{r\rightarrow 1}
b(r\gamma)$ is unique, and if $k_\gamma ^b (z) := \frac{1 - \ov{b(\gamma)}b(z)}{1-\ov{\gamma}z}$, then for all
$f \in K^2 _b$, $f(\gamma) = \ip{f}{k_\gamma ^b} _b$.
\label{thm:angderv}
\end{thm}

    With the aid of the above theorem we can now
prove the following:

\begin{thm}
    Let $V$ be a simple isometric linear transformation with deficiency indices $(1,1)$,
and characteristic function $w_V$. Any point $z \in \bm{T}$ is not an eigenvalue of any unitary extension
$U$ of $V$ if and only if the angular derivative of $w_V$ does not exist
at $z$. The symmetric linear transformation $B = \mu ^{-1} (V)$  will be densely defined if and only if
the angular derivative of $w_V$ does not exist at $z=1$. \label{thm:noone}
\end{thm}

    This proof of this theorem is almost immediately implied by the proof
of the Theorem \ref{thm:charfun} of \cite[Appendix 1]{Glazman}. We provide
a sketch of the proof here for the reader's convenience. Recall that all unitary extensions of $V$
can be labeled by a single real parameter $\alpha \in [0 , 2\pi)$,
\be U(\alpha ) := V \oplus e^{i\alpha} \ip{\cdot}{\psi _+} \psi _-\ee on $\mc{H} =
\dom {V} \oplus \dom{V} ^\perp$, where $\psi _+ \in \dom{V} ^\perp$ and $\psi _- \in \ran{V} ^\perp$
are fixed non-zero vectors of the same norm. Further recall that the vector $\psi _+$
is a generating vector for each $U(\alpha )$ \cite[Section 81, Lemma 1]{Glazman}.

\begin{proof}{ (sketch) }
    Let $w =w_V$ and let  $F(\alpha ; t ) := \chi _{[0, t)} (U (\alpha ) )$, $t \in [0 , 2\pi)$ be the spectral
distribution function of the unitary operator $U(\alpha )$. Here $\chi _\Om$ denotes the characteristic
function of the Borel set $\Om \subset [0, 2\pi]$.  Since $\psi _+$ is a generating
vector for each $U (\alpha )$, it follows that $z' = e^{i\beta} \in \bm{T}$ will be an eigenvalue for
$U(\alpha )$ if and only if the distribution function $\sigma _\alpha (t) := \ip{F(\alpha; t) \psi _+}
{\psi _+} $ has a jump at $t = \beta$. However, as in the proof of Theorem \ref{thm:charfun}, one can
calculate that
\be \frac{e^{i\alpha}}{e^{i\alpha} - w (z)} = \int _0 ^{2\pi} \frac{1}{1-e^{-it}z} d\sigma _\alpha (t),\ee
from which it follows that the value of the jump at $t=\beta$ is given by
\be \lim _{z \rightarrow e^{i\beta}} \frac{(1 - e^{-i\beta} z) e^{i\alpha}}{e^{i\alpha} - w(z)} ,\ee
where $z$ approaches $e^{i\beta }$ non-tangentially.

    In conclusion, $z' \in \bm{T}$ is not an eigenvalue of any unitary extension of $V$ if and only if
\be \lim _{z \rightarrow z'} \frac{z' -z}{e^{i\alpha} -w(z)} = 0, \ee for all $\alpha \in [0, 2\pi)$
whenever $z \rightarrow z'$ non-tangentially. This will happen if and only if
\be \lim  _{z \rightarrow z'} \left| \frac{w(z) - e^{i\alpha} } {z-z'} \right| = \infty \label{eq:dne} ,\ee for all
$\alpha \in [0 , 2\pi)$. If equation (\ref{eq:dne}) holds,
then by part (i) of Theorem \ref{thm:angderv}, the angular derivative of $w$ at $z'$ does not exist. Conversely,
suppose that the angular derivative
of $w$ does not exist at $z' = e^{i\beta}$. Observe that
\be \left| \frac{w(z) - e^{i\alpha} }{z-z'} \right| \geq \frac{1 - |w(z)|}{|z-z'|}, \ee and since
$z$ approaches $z' $ non-tangentially, $ \frac{1 - |z|}{|z' -z| }$ is bounded below in this limit. It follows
that \be  \lim  _{z \rightarrow z'} \left| \frac{w(z) - e^{i\alpha} } {z-z'} \right| \geq C \lim \inf
_{z \rightarrow z'} \frac{1 -|w(z)|}{1-|z| } = \infty, \ee  by part (iii) of Theorem \ref{thm:angderv}.\vspace{.025truein} \\

    It remains to show that $B = \mu ^{-1} (V)$ will be densely defined
if and only if $z=1$ is not an eigenvalue of any unitary extension of $V$. First, if $B$ is densely
defined, then $\ran{V-1}$ and hence $\ran{U-1}$ is dense for any unitary extension $U$ of $V$. It follows
easily from this that no such $U$ has $1$ as an eigenvalue. Conversely, assume that $B$ is not densely defined,
so that there is a vector $\xi \in \mc{H}$ such that $\ip{(V-1) \psi}{\xi} =0$ for all $\psi \in \dom{V}$. As before
choose $\psi _\pm $ such that $1 = \| \psi _\pm \| $ and $\dom{V} ^\perp = \bm{C} \{ \psi _+ \}$, $\ran{V} ^\perp =
\bm{C} \{ \psi _- \}$, and for $\ka \in \bm{C}$ define $U _\ka $ on $\mc{H} = \dom{V} \oplus \bm{C} \{ \psi _+ \}$ by
$U_\ka = V \oplus \ka  \ip{\cdot}{\psi _+} \psi _-$. Given any $ \psi \in \mc{H}$, write $\psi = \psi _V + c \psi_+$ where
$\psi _V \in \dom{V}$. Then,
\be \ip{(U_\ka -1) \psi }{\xi} = \ip{(V-1)\psi _V}{\xi}  + c \left( \ka \ip{\psi _- }{\xi} - \ip{\psi _+}{\xi} \right). \label{eq:theans}\ee

Now $\ip{\psi _-}{\xi}$ cannot be zero, as otherwise, there would exist a $\xi ' \in \dom{V}$ such that $ V \xi ' = \xi$.
This would imply $0 = \ip{(V-1) \psi }{V \xi '} = \ip{\psi}{(1-V) \xi ' }$ for all $\psi \in \dom{V}$ so that $(V-1) \xi '
= c ' \psi _+$. The fact that $\xi' \perp \psi _+$ and $\| V \xi ' \| = \| \xi ' \| $ would then imply that $V \xi ' = \xi ' $,
contradicting the simplicity of $V$.

Since $\ip{\psi _-}{\xi} \neq 0$, choose $\ka = \frac{\ip{\psi _+}{\xi}}{\ip{\psi _-}{\xi}}$ in (\ref{eq:theans}) to obtain
that $\ip{( U_\ka -1) \psi }{\xi} = 0$ for all $\psi \in \mc{H}$. This implies that $U_\ka ^* \xi = \xi$. It is straightforward
to calculate that $U _\ka ^* = V^* \oplus \ov{\ka} \ip{\cdot}{\psi _-} \psi _+$ on $\mc{H} =\ran{V} \oplus \bm{C} \{ \psi _- \}$.
Hence, writing $\xi = \xi ^* + d \psi _- $ where $\xi ^ * \in \ran{V}$ and $d \in \bm{C}$, it follows that $U_\ka ^* \xi = V^* \xi ^* + \ov{\ka} d \psi _+$,
and that \be \| \xi ^* \| ^2  + |d | ^2 = \| \xi \| ^2 = \| U_\ka ^* \xi \| ^2 = \| V^* \xi ^* \| ^2 + | \ka | ^2 |d |^2 = \| \xi ^* \| ^2
+ | \ka | ^2 | d | ^2. \ee This shows that $|\ka | =1 $ so that $U _\ka $ is the desired unitary extension of $V$ which has $1$ as an eigenvalue.

\end{proof}

    The following corollary follows readily from the above theorem, and
part (ii) of Theorem \ref{thm:angderv}.

\begin{cor}
    Let $B$ be a simple symmetric linear transformation with deficiency indices $(1,1)$,
and characteristic function $\om $.
Consider the canonical representation (\ref{eq:canon}) of $\om $. Then $B$ is densely
defined if and only if either $\sigma >0$ or \be \sum _{n \in \bm{N}} \im{z_n}
+ \intfty | t+i | ^2 d\nu (t) = \infty. \label{eq:denser} \ee If $B$ is densely defined then $\la \in \bm{R} $ is not
an eigenvalue of any self-adjoint extension of $B$ if and only if
\be \sum _{n \in \bm{N}}  \frac{\im{z_n}}{ |\la - z_n | ^2 }
+\intfty \frac{ |t +i| ^2}{|\la -t |^2} d\nu (t)  = \infty . \label{eq:noevalue} \ee
\end{cor}

The proof of this corollary is a straightforward computation.

\begin{proof}
    By Theorem \ref{thm:noone}, $B$ is densely defined if and only if the angular derivative
of $w_V$ at $\gamma =1$ does not exist. As before $V:= \mu (B)$. By part (ii) of Theorem \ref{thm:angderv} this
happens if and only if $1= \gamma \notin A_{w_V}$, \emph{i.e.} if and only if
\be \sum \frac{1 -|a_n| ^2}{| \gamma - a_n | ^2} + \int _\bm{T} | \zeta -\gamma | ^{-2} d\rho (\zeta) = \infty, \label{eq:infbaby} \ee
where $\{ a_n \}$ are the zeroes of $w_V$ in $\bm{D}$ and $\rho $ is the singular measure appearing in the
singular part of $w_V$. Explicitly, $w_V (z) = \delta \mf{B} _V (z) S_V (z)$ where $|\delta | =1$, the Blaschke part
of $w_V$ is
\be \mf{B}_V (z) = \prod _n \frac{|a_n|}{a_n} \frac{a_n -z}{1 -\ov{a_n} z} ,\ee and the singular part is
\be S_V (z) = \exp \left( -\int _\bm{T} \frac{ \zeta +z}{\zeta -z} d\rho (\zeta ) \right)
= e^{ -\frac{1+z}{1-z} \rho ( \{ 1 \} )} \exp \left( -\int _{\bm{T} \sm \{ 1 \}} \frac{ \zeta +z}{\zeta -z} d\rho (\zeta ) \right). \ee
The zeroes of $\om _B $ are $z_n = \mu ^{-1} (a_n ) \in \bm{U}$. If we let $t = \mu ^{-1} (\zeta)$ for
$\zeta \in \bm{T} \sm \{ 1\}$, $\sigma := \rho ( \{ 1 \} )$, and define the measure
$\nu $ on $\bm{R}$ by $d\nu (t) = d \rho (\mu (t)) $, then we see that the canonical representation of $\om _B$ is
$\om _B (z ) = \om _V (\mu (z) ) =  \delta \mf{B} _B (z) S _B (z)$, where the singular part $S_B$ is
\be S_B (z) = S_V (\mu (z)) = e^{i\sigma z} \exp \left( i \intfty \frac{1-tz}{t-z} d\nu (t) \right), \ee   and  $\mf{B} _B (z) = \mf{B} _V (\mu (z))$
is a Blaschke product with zero set $\{ z_n = \mu ^{-1} (a _n ) \}$.

A short calculation shows that if $\gamma =1$, then upon replacing $a_n = \mu (z_n)$,
the summand appearing in (\ref{eq:infbaby}) can be written
$\frac{1 - |a_n |^2 }{|1 -a_n |^2 } = \frac{1}{4} \left( |z_n +i |^2 - |z_n -i|^2 \right) = \im{z_n}$. Furthermore
with $t, \nu $ and $\sigma$ as above, and $\gamma =1$  the integral appearing in (\ref{eq:infbaby}) can be written
\be \int _\bm{T} | \zeta - 1| ^{-2} d\rho (\zeta) = + \infty \cdot \sigma + \int _{\bm{T} \sm \{1 \} } | \zeta -1 | ^{-2} d\rho (\zeta)
= +\infty \cdot  \sigma + \intfty |t+i| ^2 d \nu (t) .\ee Here $+\infty \cdot \sigma $ is defined as $+\infty$
if the point mass $\sigma = \rho (\{ 1 \}) >0$ and $0$ if $\sigma =0$.
This proves the statement characterizing when $B$ is densely defined.

The second statement follows from similar calculations, and the observation that $\la \in \bm{R}$ is an eigenvalue of
a self-adjoint extension of $B$ if and only if $\gamma = \mu (\la) \in \bm{T} \sm \{ 1 \}$ is an eigenvalue of some unitary
extension of $V= \mu (B)$.
\end{proof}

\begin{cor}
    The operator, $M ^\phi$ of multiplication by $z$ in $K^2 _\phi$ is densely defined if and only
if the angular derivative of $\varphi =\phi \circ \mu ^{-1}$ does not exist at $z=1$. If $M ^\phi$ is densely defined,
then $\la \in \bm{R}$ is not an eigenvalue of any of its self-adjoint extensions if and only if
the angular derivative of $\varphi$ at $\mu (\la)$ does not exist. \label{cor:densemult}
\end{cor}

\begin{proof}
    Using the formula for $w _\varphi$ given in equation (\ref{eq:charfunction}), it is not difficult
to verify that the angular derivative of $w_\varphi$ at $z \in \bm{T}$ exists if and only if the angular derivative
of $\varphi$ exists at $z$.
\end{proof}

\section{Spectra of self-adjoint extensions of $M^\phi$}

\label{section:spectra}

    In this section we consider the case where $\phi$ is an extreme point, and
$\phi$ obeys the conditions of Corollary \ref{cor:densemult} so that $M ^\phi$, which
acts as multiplication by $z$, is a closed, simple and densely defined symmetric linear operator with
deficiency indices $(1,1)$ in $K^2 _\phi$.
The following result, due to Lifschitz, immediately identifies the essential spectrum
of $M^\phi$ in $K^2 _\phi$. Recall that a point
$z \in \bm{C} $ is said to be a regular point, or a point of regular type for a closed linear transformation $T$ if
$T-z$ is bounded below on $\dom{T}$. A symmetric linear transformation $B$ is said to be regular
if every $z \in \bm{C}$ is regular for $B$.

\begin{thm}{ (Lifschitz) }
    In order that a real number $\la $ be a point of regular type of a simple
symmetric linear transformation $B$ with deficiency indices $(1,1)$ and characteristic function $\om$,
it is necessary and sufficient that both of the following conditions be satisfied:
\bn
    \item  $\om$ is analytic in a neighbourhood of $\la $.
    \item $| \om (t) | =1 $ on some open interval of $\bm{R} $ containing
    $\la$.
\en
\label{thm:regular}
\end{thm}

In \cite[Theorem 4]{Lifschitz2}, a more general version of the above theorem is established for simple isometric
operators with deficiency indices $(n,n)$, $n < \infty$.

\subsubsection{Essential spectrum of $M ^\phi$}

\label{subsubsection:esspec}

For a simple symmetric operator with deficiency indices $(1,1)$,
$\la $ belongs to the essential spectrum of $B$ if and only if $\la $
is not a regular point of $B$. For the operator of multiplication $M ^\phi$,
its characteristic function
$\om _\phi $ obeys the conditions of Theorem \ref{thm:regular} at $\la \in \bm{R}$
if and only if $\phi $ obeys those same conditions. It follows that
$\sigma _e (M ^\phi) = \supp{\phi} \cap \bm {R} \cup \{ \infty \} $
where $\supp{\phi}$ is
defined as the union of the closure of the set of zeroes of $\phi$ and
the closed support of the measure $\nu$ that appears in the canonical
representation of $\phi$ in equation (\ref{eq:canon}). This is clear as if
$\la \in \bm{R}$ is either a limit point of the zeroes of $\phi$ or in the
closed support of the part of the measure $\nu$ which is singular with respect
to Lebesgue measure, then $\phi$
does not satisfy condition (1) of Theorem \ref{thm:regular} at $\la$. Now suppose
that $\la $ belongs to the closed support of the absolutely continuous part of
$\nu$. By the inner-outer factorization for the $H^\infty $ function $\phi$,
the absolutely continuous part of $\nu$ is given by $d\nu (\gamma)  = -\ln | \phi (\gamma ) | dm (\gamma) $, where
$m$ is normalized Lebesgue measure on $\bm{T}$.  Since  $\phi \in B _1 (H^\infty(\bm{U}))$, it follows that
a Borel set $\Om \subset \bm{T}$ belongs to the closed support of the absolutely continuous part of $\nu$
if and only if $| \phi ( \gamma ) | < 1 $ almost everywhere for $\gamma \in \Om$.   It follows that $\la$
is the limit of a sequence $\la _n \in \bm{R}$ where $| \phi (\la _n ) | < 1 $   for each $n$,
so that each $\la _n$ belongs to $\sigma _e (M ^\phi)$ by condition (2) of Theorem
\ref{thm:regular}. Since the essential spectrum is closed, $\la \in \sigma _e (M ^\phi)$
as well.

    If $A'$ is any self-adjoint extension of a symmetric operator $A$ with
finite deficiency indices, then $\sigma _e (A') = \sigma _e (A)$. This follows easily
from the fact that $\mu (A')$ is a finite rank extension of $\mu (A)$.

\subsubsection{Remark} Note that since $M ^\phi$ is simple, it has no eigenvalues so that $\sigma (M ^\phi) = \sigma _e (M ^\phi)$.
Further observe that if $\la \in \bm{R}\sm\sigma _e (M ^\phi)$, then by Theorem \ref{thm:regular},
the angular derivative of $\om _\phi$ (and hence of $\phi$) exists at $\la$. If   \label{subsubsection:ang}

\subsubsection{Total orthogonal sets of point evaluation vectors} The spectra of self-adjoint
extensions of $M ^\phi$ are related to the existence of total orthogonal sets of point evaluation vectors
in $K^2 _\phi$. As discussed in \cite{Fricain}, a set of point evaluation vectors $\Gamma:= \{ k ^\phi _{\la _n} \} _{n \in \bm{Z}}$
in $K^2 _\phi$ can only be orthogonal if $\{ \la _n  \} \subset \bm{R}$. If $\Gamma $ is a total
orthogonal set, the set of points $ \{ \la _n \}$ is the spectrum of a self-adjoint
extension $M'$ of $M ^\phi$, and  $\Gamma$ consists of eigenvectors to $M'$. To see this note
that since $M^\phi$ is densely defined, each $k_{\la _n} ^\phi$ is an eigenvector of $(M^\phi)^*$
with eigenvalue $\la _n$. Hence the closure, $M'$ of the restriction of $(M^\phi )^*$ to the linear
span of $\Gamma$ is a self -adjoint restriction of $(M^\phi)^*$. If $\psi \in \dom{M}$, then
$\ip{\psi}{M' \psi '} = \ip{\psi}{(M^\phi) ^* \psi ' } = \ip{M^\phi \psi}{\psi '}$ for all $\psi ' \in
\dom{M '}$. This implies that $\psi \in \dom{ (M') ^* = M'}$ and $(M')^* \psi = M' \psi = M ^\phi \psi$,
so that $\dom{M ^\phi} \subset \dom{M'}$, and $M'$ is a self-adjoint extension of $M^\phi$.

    In particular, Theorem \ref{thm:regular} and Remark \ref{subsubsection:ang} allow one to conclude that if
$K^2 _\phi$ has a total orthogonal set of point evaluation vectors, then the angular derivatives
of $\varphi$ must exist almost everywhere on $\bm{T}$ so that $\phi $ must be inner. This follows as if the angular
derivative of $\varphi $ does not exist at a point $z \in \bm{T}$, then either $\varphi$ does not have
unit modulus at $z$ or it is not analytic in any neighbourhood of $z$ so that $\mu ^{-1} (z)$ belongs
to the essential spectrum of $M ^\phi$. If there is a Borel subset $\mu (\Om ) \subset \bm{T}$ of non-zero measure on which
the angular derivatives of $\varphi $ do not exist, then $\Om := \mu ^{-1} (\mu (\Om) ) \subset \bm{R}$ belongs
to $\sigma _e (M ^\phi)$ and no point of $\Om$ belongs to the point spectrum of any self-adjoint extension of
$M ^\phi$ by Corollary \ref{cor:densemult}. Hence if $M' $ is any self-adjoint extension of $M ^\phi$, and $\La := \{ k _{\la _n} \}$ is
a maximal set of orthogonal eigenvectors to $M'$, then $\Om \subset \sigma _e (M) = \sigma _e (M')$ and
$\chi _\Om (M')$ projects onto a non-zero subspace
orthogonal to the closed linear span of $\La$. Hence $K^2 _\phi$ has no total orthogonal set of point
evaluation vectors. This fact appears as Corollary 2.2 in \cite{Fricain}.

\label{subsubsection:orthog}
\subsection{Spectra of self-adjoint extensions of $M^\phi$}
\label{subsection:pspec}

    Suppose that $\la \in  \bm{R}\sm\sigma _e (M^\phi)$. By Remark \ref{subsubsection:ang},
the angular derivative of $\phi$ at $\mu (\la)$ exists, and it follows from Theorem \ref{thm:angderv} that
point evaluation at $\la$ is a bounded linear functional on $K^2 _\phi$,
generated by the point evaluation vector
\be k_\la ^\phi := \frac{i}{2\pi} \frac{1 -\ov{\phi (\la)} \phi }{z -\la}. \ee

    The spectra of any fixed self-adjoint extension, $M'$ of $M^\phi$ is
$\sigma (M') = \sigma _p (M') \cup \sigma _e (M ^\phi)$ where $\sigma _p (M')$
denotes the set of eigenvalues of $M'$. Note that if $\la \in \bm{R} \sm \sigma _e (M ^\phi)$
is an eigenvalue of $M'$, then it follows that $k_\la ^\phi$ is an eigenvector of
$M'$ to eigenvalue $\la$. To see this, first note that since $M ^\phi$ is densely
defined, $k_\la ^\phi$ will be an eigenvector of $(M^\phi) ^*$ to eigenvalue $\la$.  If $\la$ is an
eigenvalue of $M'$ and if $f _\la $ is an eigenvector
for $M'$ to eigenvalue $\la$, it must be that $f_\la = c k ^\phi _\la $ for some $c \in \bm{C}$ since $\ker{(M^\phi) ^* -\la}
= \bm{C} \{ k ^\phi _\la \}$ is one dimensional and $(M^\phi) ^*$ extends $M'$.

   Given the Cayley transform $V  ^\phi  =\mu (M^\phi) = (M^\phi-i) (M^\phi +i )$ of $M^\phi$, we have that
$\dom{V ^\phi} = \ran{M^\phi+i}$ and $\ran{V _\phi } = \ran{M^\phi-i}$. Let $\psi _{\pm }$ be fixed non-zero vectors
in $\ran{M^\phi \mp i } ^\perp = \ker{(M^\phi) ^* \pm i } $ which have the same norm, $\| \psi _+ \|
= \| \psi _{-} \| $. Recall that the family of all unitary extensions of $V ^\phi$ can be
labeled by a single real parameter $\alpha \in [0 , 2\pi)$ as follows (see, for example \cite{Kempfsamp},
\cite{Martin-symsamp}, or \cite{Martin-uncer}). Given any such
an $\alpha $, define \be U ^\phi (\alpha ) := V ^\phi \oplus e^{i\alpha} \ip{\cdot}{\psi _+ } _\phi \psi _- . \label{eq:unix} \ee
As $\alpha $ ranges in $[0, 2\pi)$, $U ^\phi (\alpha )$ covers all possible unitary extensions of
$V ^\phi$, and the set of all $M ^\phi (\alpha ) := \mu ^{-1} ( U ^\phi  (\alpha ) )$ for $\alpha \in [0, 2\pi)$
is the family of all self-adjoint extensions of $M ^\phi$. In what follows, we choose $\psi _{-}= i 2\pi k _i ^\phi$
and $\psi _+ = - C_\phi \psi _{-}$. Replacing $\psi _\pm$ by $c \psi _\pm$ where $c \neq 0$ does not change $U(\alpha)$.

    The domain of $M ^\phi (\alpha )$ can then be decomposed as:
\be \dom{M ^\phi (\alpha )} = \dom{M ^\phi } \dotplus \bm{C} \left\{ \psi _+  + e^{i\alpha} \psi _{-} \right\} \ee
\cite[Section 80]{Glazman}, where $\dotplus $ denotes the non-orthogonal direct sum of
linearly independent subspaces, and $\bm{C} \{ \psi \}$ is the one dimensional linear span of a vector
$\psi$. It follows that a point $x \in \bm{R} \sm \sigma _e (M^\phi)$ will belong to $\sigma _p (M^\phi (\alpha))$
if and only if there is a non-zero $c \in \bm{C}$ such that $f = -i 2\pi k_x - c ( \psi _+ + e^{i\alpha} \psi _{-} )
$ belongs to $\dom{M ^\phi}$. If $f \in \dom{M^\phi }$ then for any $z \in \bm{U}$ it follows that
$\ip{M^\phi f}{k_z ^\phi }_\phi =z f(z)$. Alternatively, \be M^\phi f (z)  = (M^\phi) ^*f (z) = -i 2\pi x k_x (z)
 -ic ( \psi _+ (z) - e^{i\alpha} \psi _{-} (z) ). \ee
Equating these two expression for $M^\phi f (z)$ yields
\ba & &  z \frac{1 -\ov{\phi (x)} \phi (z)}{z-x}  - z c  \frac{\phi (z) + \phi (i) } {z-i} - zc e^{i\alpha}
\frac{1 - \ov{\phi (i)} \phi (z) }{z+i} \nonumber \\ & = & x \frac{1 -\ov{\phi (x)} \phi (z)}{z-x} -i
c  \frac{\phi (z) - \phi (i) } {z-i} -i c e^{i\alpha}
\frac{1 - \ov{\phi (i)} \phi (z) }{z+i}. \ea
    This expression can be simplified to yield
\be 0 = \left( 1 + c \phi (i) +c e^{i\alpha} \right) - \phi (z) \left( \ov{\phi (x)} +c -c \ov{\phi (i)} e^{i\alpha} \right).\ee
Since this must hold for all $z \in \bm{U}$, both of the bracketed terms must vanish separately, leading to the expression
\be \phi (x) = \frac{ e^{i\alpha} - \phi (i) }{1 - \ov{\phi (i) } e^{i\alpha}} \label{eq:pspec}.\ee Note that the right hand
side of this formula has modulus $1$, as it must since the angular derivative of $\varphi $ at $\mu (x)$ exists.

\subsubsection{Summary of Results}   In summary, we have proven that $\sigma (M ^\phi (\alpha) ) = \sigma _e (M^\phi ) \cup \sigma _p (M^\phi (\alpha ))$
where $ \sigma _e (M^\phi ) = \supp{\phi} \cap \bm{R} $ and $\sigma _p (M^\phi (\alpha )) =
\{ x \in \bm{R} \sm \sigma _e (M^\phi )  | \  \phi (x)  = \frac{ e^{i\alpha} - \phi (i) }{1 - \ov{\phi (i) } e^{i\alpha}} \}.$ In
particular if $\phi (i) =0$ it follows that $\sigma _p (M ^\phi  (\alpha ) ) = \{ x \in \bm{R}\sm
\sigma _e (M^\phi ) | \ \phi (x) = e^{i\alpha} \}$. Here recall that $M^\phi (\alpha )$ is the inverse Cayley transform of $U ^ \phi (
\alpha )$, as given in (\ref{eq:unix}), with the specific choice of deficiency vectors $\psi _- = k_i ^\phi$, $\psi _+ = - C_\phi \psi _-$.
\label{subsubsection:sum}

\subsubsection{Remark} It is a simple calculation to verify that if $x,y \notin
\sigma _e (M)$, $x \neq y$, then $\ip{k_x ^\phi}{k_y ^\phi} _\phi =0$ if and only if $\phi (x) = \phi (y)$. Indeed,
this inner product is equal to $k_x ^\phi (y) = \frac{1}{2\pi} \frac{1-\ov{\phi (x)} \phi (y) }{y-x}$,
and since $x,y \notin \sigma _e (M)$, $|\phi (x) | = | \phi (y) | =1$, by Remark \ref{subsubsection:ang}.

\section{de Branges spaces}
\label{section:dB}

    Recall that a de Branges function $E$ is an entire function which obeys
$|E(z) | > | E(\ov{z} )| $ for all $z \in \bm{U}$. Given such a function the de
Branges space $\mc{H} (E)$ is defined as the space of all entire functions $f$
such that $f/E$ and $f^*/E$ belong to $H^2 (\bm{U})$ \cite[Prop. 2.1]{Remling}. The space $\mc{H} (E)$ is a
Hilbert space with respect to the inner product $\ip{f}{g} _{\mc{H} (E)} =
\ip{f/E}{g/E} _{L^2 (\bm{R} )}$ \cite{deBranges}.

    Suppose that $\phi = F $ is an inner function which can be extended to a meromorphic function
on $\bm{C}$. An inner function $F$ has this form if and only if its singular part has support only at the point
at infinity, \emph{i.e.} the singular part consists only an exponential term $e^{i\sigma z}$, $\sigma \geq 0$,
and the zeroes of $F$ have no finite accumulation point on $\bm{R}$. Recall that an inner function
$F$ has these properties if and only if there exists
a de Branges function $E$ such that $F = \frac{E^*}{E}$ (Theorem \ref{thm:innerdB}). The following results characterizing
de Branges spaces and de Branges functions can be found, for example, in \cite{Baranov} \cite{Havin}.

\begin{lemming}
The map $U : \mc{H} (E) \rightarrow K^2 _{E^*/E}$ defined by $U f  = f / E$ is an isometry from $\mc{H} (E) $ onto $K^2 _{E^*/E}$.
\label{lemming:isodB}
\end{lemming}

\subsubsection{Remark} It is clear that if $E$ is a de Branges function, then $E^*/E$ is a meromorphic inner function. The defining
inequality $|E(z)| > |E (\ov{z} )| $ for $z \in \bm{U}$ ensures that the zeroes of $E$ are contained in the closure of the lower half
plane and that $E^*/E$ is bounded by $1$ in $\bm{U}$, and is unimodular on $\bm{R}$. Since $E$ is entire, its zeroes have no finite accumulation
point, so that $E^*/E$ is analytic on some neighbourhood of each $x \in \bm{R}$, and meromorphic in $\bm{C}$.

    If $w \in \bm{L}$ then $z-w$ is a de Branges function, also $e^{-i\sigma z}$ is a de Branges function for any $\sigma >0$. Any
finite product of de Branges functions, or of a de Branges function with any entire function $G$ such that $G=G^*$, and
$G$ has only real zeroes, is again a de Branges function. It follows easily that if $F(z)= e^{i\sigma z} \mf{B}(z)$ is a meromorphic inner function where $\mf{B}(z) = \prod _{n=1} ^N \frac{ z - z_n}{z- \ov{z_n}} $ is a finite Blaschke product, then
$E(z) = \gamma G(z) e^{-i\sigma z /2}  \prod _{n=1} ^N (z -\ov{z_n})$, is a de Branges function satisfying $F =E^*/E$.
Here, $\{ z_n \} \subset \bm{U}$, $G$ is as described previously, and $\gamma \in \bm{T}$. The following theorems
generalize these results to the case where $\mf{B}$ is an infinite product.

\begin{thm}{ (M.G. Krein) }
    If $E$ is a de Branges function and $(\ov{z} _n ) _{n\in \bm{N}}$ are its zeroes in $\bm{L}$ ordered so that $|z_n | \leq | z_{n+1} |$, then $\sum _{n=1} ^\infty \left|
\im{\frac{1}{ z_n }} \right| < \infty$ and,
\be E(z) = \gamma G(z) e^{-i\sigma z } \prod _{n=1} ^\infty \left( 1 - \frac{z}{ \ov{z_n} } \right) e^{\frac{1}{2} ( p_n (z) + p_n ^* (z) )}, \label{eq:dBform}
 \ee
where $\sigma >0$, $| \gamma | =1$, $p_n (z) := \sum _{k=1} ^n \frac{1}{k \ov{z_n} ^k} z^k$ and $G=G^*$ is an entire function whose zeroes lie on the real axis. \label{thm:dBfunrep}
\end{thm}

\begin{thm}
    If $F \in H^\infty (\bm{U})$ is inner, then $F=E^*/E$ for some de Branges function $E$
if and only if $F(z) = e^{i\sigma z} \mf{B} (z)$ ; $\sigma \geq 0$,
where $\mf{B} (z)$ is a Blaschke product whose zeroes have no finite accumulation point.

Let $F=e^{i\sigma' z} \mf{B} (z)$ be a meromorphic inner function with zeroes
$\{ z_n \} _{n=1} ^\infty $, and $\sigma' \geq 0$. A de Branges function $E$
satisfies $F = E^* /E $ if and only if both $\{ \ov{z_n} \} _{n=1} ^\infty$ is the
set of non-real zeroes of $E$, and the constant $\sigma$ appearing in the canonical
representation (\ref{eq:dBform}) of $E$ is equal to $\sigma ' /2$. Hence, such an $E$
is determined uniquely by $F$ up to a unimodular constant $\gamma$ and an entire function $G$ which obeys $G=G^*$,
and whose zeroes lie on $\bm{R}$.

\label{thm:innerdB}
\end{thm}

    These facts lead to the following representation theorem
for regular simple symmetric operators with deficiency indices $(1,1)$.

\begin{thm}
    Let $B$ be a linear transformation with domain and range contained in a separable
Hilbert space $\mc{H}$. Then $B$ is regular, closed and simple symmetric with deficiency deficiency indices $(1,1)$
if and only if it is unitarily equivalent to multiplication by $z$ in a de Branges space of entire functions. Equivalently,
such a $B$ is unitarily equivalent to multiplication by $z$ in $K^2 _F$ where $F \in H^\infty (\bm{U})$ is a meromorphic
inner function.
\label{thm:dBrep}
\end{thm}

    By multiplication by $z$ in a de Branges space $\mc{H} (E)$, we mean the linear
transformation in $\mc{H} (E)$ which acts as multiplication by $z$ on its domain,
and which has no proper extension. The fact that multiplication by $z$ in any de Branges space $\mc{H} (E)$ is a closed regular simple
symmetric linear transformation with deficiency indices $(1,1)$ is well known, see for example
\cite[Theorems 16-17]{Martin-symsamp} \cite{deBranges}. It remains to prove necessity.

\begin{proof}
    Suppose $B$ is closed, regular and simple symmetric with deficiency indices $(1,1)$.
If $\om _B$ is the characteristic function of $B$, then the regularity of $B$ and
Theorem \ref{thm:regular} imply that $|\om _B (\la ) | =1$, and that $\om _B$ is analytic in
a neighbourhood of $\la$ for any $\la \in \bm{R}$. It follows
that $\om _B$ is a meromorphic inner function and has the form $\om _B (z) = e^{i\sigma z} \mf{B} (z)$ where $\sigma \geq 0$,
and $\mf{B} (z)$ is a Blaschke product whose zeroes have no finite accumulation point. By
Theorem \ref{thm:innerdB}, there is a de Branges function $E$ such that $\om _B = E^*/E$.
By Theorem \ref{thm:symrep}, $B$ is unitarily equivalent to multiplication by $z$ in $K^2 _{\om _B }$, and
since by Lemma \ref{lemming:isodB}, multiplication by $E$ is an isometry from $K^2 _{\om _B }$
onto the de Branges space $\mc{H} (E)$, it follows that $B$ is unitarily equivalent to
multiplication by $z$ in $\mc{H} (E)$. The converse follows from the comment preceding this proof.
\end{proof}

    More precisely, by \ref{thm:symrep}, the following is true.

\begin{cor}
    If $B$ is a regular simple symmetric linear transformation with deficiency indices $(1,1)$,
then its characteristic function $\om _B$ is a meromorphic inner function. For $\alpha \in \bm{D}$, let
$\om _B ^{(\alpha)} :=
\frac{\om _B -\alpha}{1- \ov{\alpha} \om _B}$. If $e^{i\frac{\sigma (\alpha)}{2} z }$ is
the singular part of $\om _B ^{(\alpha)}$ and $(z_n (\alpha) ) _{n=1} ^\infty$ are its zeroes, then $B$ is
unitarily equivalent to multiplication by $z$ in
any de Branges space $\mc{H} (E _\alpha )$ where
\be     E _\alpha (z) = \gamma (\alpha) G _\alpha (z) e^{-i\sigma (\alpha) z } \prod _{n=1} ^\infty \left( 1 - \frac{z}{ \ov{z_n (\alpha)} } \right) e^{\frac{1}{2} ( p_n (\alpha ; z) + p_n ^* (\alpha ; z) )}, \ee
$\sigma (\alpha) >0$, $| \gamma (\alpha) | =1 $, $p_n (\alpha ; z) := \sum _{k=1} ^n \frac{1}{k \ov{z_n (\alpha) } ^k} z^k$ and $G _\alpha$ is any entire function whose zeroes lie on the real axis and which obeys $G _\alpha =G _\alpha ^*$.
\end{cor}

The results of Section \ref{section:dense} applied to this particular case when
$\phi = F = E^*/E $ for some de Branges function $E$ yields the following
criterion for multiplication by $z$ to be a densely defined symmetric operator
in $\mc{H} (E)$.

\begin{thm}
    Let $E$ be a de Branges function, and consider the representation of $E$
given in Theorem \ref{thm:dBfunrep}. Then multiplication
by $z$ in $\mc{H} (E)$ is densely defined if and only if at least one
of the following two conditions holds:
(i) $\sigma > 0$, or
(ii) $\sum _{n \in \bm{N}} \im{z_n} = \infty$.

\label{thm:dBdense}
\end{thm}

\subsection{The spectra of self-adjoint extensions of $M ^F$}

    Since the operator of multiplication, $M$, in $K^2 _F$, where $F$ is
a meromorphic inner function,
has no essential spectrum, the spectra of any of its self-adjoint extensions
is purely discrete. Recall that any symmetric linear transformation $B$ for which $B-z$ is bounded below
for all $z \in \bm{C}$ is said to be regular. Since $M ^F$, where $F$ is a meromorphic
inner function is simple and $\sigma _e (M ^F) = \{ \infty \}$, it follows that $M^F$
is regular.  Let the self-adjoint extension $M^F (\alpha)$ of $M^F$ be defined as in
Subsection \ref{subsection:pspec}, \emph{i.e.} choose the deficiency vectors $\psi _- = k_i ^F$
and $\psi _+ = -C_F \psi _-$.  Since every extension
has an infinite number of eigenvalues, it follows from equation (\ref{eq:pspec})
of Subsection \ref{subsection:pspec} for the spectrum of $M ^F (\alpha )$ that given any $\beta
\in [0 ,2\pi)$ there is an infinite number of points $x \in \bm{R}$ such that $F(x) = e^{i\beta}$.

\subsubsection{Definition} \label{subsubsection:bdef} For each $\beta \in [0, 2\pi)$, choose a point $x \in \bm{R}$ such that $F(x) = e^{i\beta}$
and let $M ^F _\beta $ be that self-adjoint extension of $M ^F$ for which $x \in \sigma (M^F _\beta )$. \vspace{.05truein} \\

The results of Section \ref{section:spectra} show that $M^F _\beta$ is well defined, and that
 \be \sigma (M ^F _\beta ) = \{ x \in \bm{R} |
\ F(x) = e^{i\beta} \}. \label{eq:betaspec} \ee Furthermore by inverting equation (\ref{eq:pspec}) of Section \ref{section:spectra},
$M ^F _\beta = M^F (\alpha )$ where  $e^{i\alpha} = \frac{ F (i) + e^{i\beta} }{1 + \ov{F(i)} e^{i\beta}} $.
In particular, if $F(i) = 0 $, $M^F _\beta = M  ^F (  \beta )$.

    Since $F$ is unimodular on $\bm{R}$, we have $F(x) = e^{i\tau (x)} = E^* (x) / E(x) $ for a real-valued
function $\tau$. It follows that
the function $\tau : \bm{R} \rightarrow \bm{R} $ can be defined so that it is infinitely differentiable and has a local
analytic extension about any point $x \in \bm{R}$. We will call such a function $\tau$ a phase function of $F$.
Observe that the spectrum of $M ^F _\beta$ can
be written \be \sigma (M^F _\beta ) = \{ x \in \bm{R} | \  \tau (x) = \beta + 2\pi n ; \ n \in \bm{Z} \}.
\label{eq:betaspec2} \ee This
implies in particular that any $x \in \bm{R}$ belongs to the spectrum of some self-adjoint extension of $M^F$.
More precisely, the following results hold :

\begin{thm}{ (\cite[Theorem 2]{Martin-symsamp}, \cite{Kempfsamp}) }
    Let $B$ be a closed symmetric operator densely defined in
$\mc{H}$. If $B$ is simple, regular and has deficiency indices
$(1,1)$, then the spectra of any one of its self-adjoint
extensions consists of eigenvalues of multiplicity one with no
finite accumulation point. Furthermore, the spectra of all of its
self-adjoint extensions cover $\bm{R}$ exactly once.
\label{thm:once}
\end{thm}

\begin{thm}
    If $\tau :\bm{R} \rightarrow \bm{R}$ is such that $F(x) = e^{i\tau (x)}$,
then $\tau ' (x) = 2 \pi \| k_x ^F \| ^2   > 0$ and $\tau $ is a $C^\infty$ bijection
of $\bm{R} $ onto $(-\infty , b)$, $(a, \infty)$ or $\bm{R}$, depending on whether
the spectrum of each self-adjoint extension of $M^F$ is bounded above, bounded below
or neither bounded above nor below, respectively. The phase function $\tau$ has a local analytic extension about any point $x \in \bm{R}$.
\label{thm:inv}
\end{thm}

    Note that if $E$ is a de Branges function such that $F=E^*/E$ then $\tau (x)/2$ is a phase function associated
with $E$, as defined in \cite{deBranges}. The fact that if one self-adjoint extension of $M^F$ is bounded above or below, then all are follows
immediately from Krein's alternating eigenvalue theorem (\cite{Krein}, pg. 19):

\begin{thm}{ (Krein) }
    Let $B$ be a closed simple symmetric operator in $\mc{H}$
with deficiency indices $(1,1)$. Suppose that the interval $I
\subset \bm{R}$ consists of regular points of $B$. Then, the
eigenvalues of any two self-adjoint extensions $B ' $ and $B''$ of
$B$ alternate in $I$. \label{thm:alternate}
\end{thm}

\begin{proof}{ (of Theorem \ref{thm:inv}) }
    Since $F$ is analytic on a region containing $\bm{R}$, Theorem \ref{thm:angderv} implies that
point evaluation at any point $x \in \bm{R}$ is
a bounded linear functional in $K^2 _F$; given any $f \in K^2 _F$, and any $x \in \bm{R}$, $\ip{f}{k_x} = f(x)$,
where $k_x ^F  (z) = \frac{i}{2\pi} \frac{1 - \ov{F(x)}F(z)}{z-x}$. It is straightforward to calculate
that $ 0 \leq \| k _x ^F \| ^2 = k_x (x) = \frac{1}{2\pi i} \ov{F (x)} F'(x) $, and since $e^{i\tau (x) } = F(x)$
it follows that $\tau ' (x) = -i \ov{F(x)} F'(x) = 2\pi  \| k_x ^F \| ^2 \geq 0$.  To show that
$\tau $ is strictly increasing, and hence injective, it remains to show that $\| k _x ^F \| >0$ is strictly
positive for every $x \in \bm{R}$. To see that $k_x ^F \neq 0$ for any $x \in \bm{R}$, recall that by definition,
$k_x ^F  (z) = \frac{i}{2\pi} \frac{1 - \ov{F(x)}F(z)}{z-x}$ which is non-zero almost everywhere
with respect to Lebesgue measure since $F$ is a non-constant inner function. Hence $\| k _x \| >0$,
and $k_x \neq 0$ for any $x \in \bm{R}$.

    By (\ref{eq:betaspec2}) and the fact that $\tau $ is strictly increasing on $\bm{R}$,
the spectra of each self-adjoint extension is bounded
above or below if and only if the range of $\tau$ is bounded above or below. The phase
function $\tau$ cannot be bounded both above and below as this, and the fact that
the spectra of each self-adjoint extension of $B$ has no finite limit point (by Theorem \ref{thm:once})
would imply that each
self-adjoint extension has only a finite number of eigenvalues. Since each such self-adjoint
extension is unbounded this is not possible.
Now suppose that $\tau $ is bounded above and that $b = \sup _{x \in \bm{R}} \tau (x) =
\lim _{x \rightarrow \infty} \tau (x)$. Then since $\tau $ is strictly increasing,
$b$ is not in the range of $\tau$ and $\tau $ is onto $(-\infty , b)$. The other two
cases are similarly easy to verify.
\end{proof}

\subsubsection{Remark} The spectrum of $M ^F _\beta$ is $\sigma (M^F _\beta) = \{ x \in \bm{R} |
\ \tau (x) = \beta + 2\pi n ; \ n \in \bm{Z} \cap \ran{\tau} \} $.
Each $M^F $ is unitarily equivalent to $M^{\gamma F}$ where $\gamma$
is any unimodular constant. If $\tau _\gamma $ is a phase function of $\gamma F$, $\gamma F (x) = e^{i\tau _\gamma (x)}$,
then there is always a $\gamma \in \bm{T}$ such that $\ran{\tau _\gamma} = (-\infty , 0 )$, $(0, \infty )$ or $\bm{R}$, and $\ran{\tau  _\gamma} \cap
\bm{Z} = \pm \bm{N}$ or $\bm{Z}$. For example, if $\ran{\tau} = (a, \infty )$, let $\gamma = e^{-ia}$.
Then $\tau _\gamma = \tau -a $ is a phase function for $\gamma F$ with range $(0, \infty)$. \label{subsubsection:range}

\begin{cor}
    Given a phase function $\tau$ for $F$, let
$\la : \bm{R} \rightarrow \bm{R}$ denote the monotonically strictly increasing $C^\infty$
function of $\ran{\tau}$ onto $\bm{R}$ which is the inverse of $\tau$, $\la (\tau (x) ) = x$.
Then $\sigma (M _\beta ) = ( \la (\beta + 2\pi n ) ) _{n \in \bm{Z} \cap \ran{\tau} }$.
\end{cor}

    Given an arbitrary densely defined simple symmetric operator $B$ with deficiency indices
$(1,1)$, recall that, as in Subsection \ref{subsection:pspec}, all self-adjoint extensions $B(\alpha)$ of $B$
can be labeled by a single parameter $\alpha \in [0 ,2\pi)$ where $B(\alpha) := \mu ^{-1} (U(\alpha ) )$ is defined as the
inverse Cayley transform of \be U(\alpha ) = \mu (B) \oplus e^{i\alpha} \ip{\cdot}{\psi _+}\psi _-,\ee
and $\psi _\pm \in \ker{B^* \pm i} = \bm{C} \{ \psi _\pm \}$ are chosen so that $\| \psi _+ \| = \| \psi _- \| \neq 0$.

\begin{cor}
    Let $B$ be a regular simple symmetric linear operator with deficiency indices $(1,1)$ and
characteristic function $\om$. Let $\tau$ be a phase function for $\om$. Then $\tau$ is a $C^\infty $ strictly monotonically increasing diffeomorphism of $\bm{R} $ onto
its range, $\tau ' (x) >0 $ for all $x \in \bm{R}$, and $\tau$ has
a local analytic extension about any point $x \in \bm{R}$. If $B _\beta$, $\beta \in [0,2\pi )$ is defined as that self-adjoint
extension of $B$ such that there is an $x \in \sigma (B_\beta )$ so that
$\om (x) = e^{i\beta }$, then $\sigma (B _\beta )
= \{ \la (\beta + 2\pi n ) \} _{n \in \bm{Z} \cap \ran{\tau} }$, where $\la = \tau ^{-1}$. If $B(\alpha )$ is
defined as above, then there is a $\theta \in [0 ,2\pi)$ such that $\beta = \alpha +\theta \ \mr{mod} \ 2\pi$,
and $\psi _\pm $ can be chosen so that $\beta =\alpha$, \emph{i.e.} so that $B _\alpha = B (\alpha)$. \label{cor:specfunk}
\end{cor}

The above corollary shows that the spectra of the self-adjoint extensions
of $B$ behave very smoothly with respect to the parameter $\alpha$ labeling the self-adjoint extensions.

\begin{proof}
    The bulk of the corollary follows immediately from the previous results of this section, and the fact that $B$
is unitarily equivalent to $M ^\om$ in $K^2 _\om$. Assume, without loss of generality that $B = M ^\om$.
We will verify the final assertion. By the results \ref{subsubsection:sum} of
Subsection \ref{subsection:pspec}, and Definition \ref{subsubsection:bdef}, with the choice $\psi _- = k_i ^\phi$
and $\psi _+ = -C_\phi \psi _-$, $M^ \om _\beta = M ^\om (\beta )$ for all $\beta \in [0, 2 \pi)$. If one makes
a different choice $\psi _\pm ' $ of deficiency vectors then $\psi _\pm ' = c \chi _\pm \psi _\pm $, where $\psi _\pm$ are
as before, $c = \frac{ \| \psi _\pm ' \| }{ \| \psi _\pm \| } \neq 0$ and  $\chi _\pm \in \bm{T}$. It follows that with this
choice of deficiency vectors, if $\ov{\chi} _+ \chi _- = e^{i\theta}$, then \be U ^\om (\alpha ) = V \oplus e^{i\alpha} \ip{\cdot}{\chi _+ \psi _+ }
\chi _- \psi _- = V \oplus e^{\alpha + \theta} \ip{\cdot}{\psi _+}{\psi _-} = U ^\om (\beta),\ee
where $U^\om (\alpha ) = \mu (M^\om (\alpha) )$ and $\beta = \alpha + \theta \ \mr{mod} \ 2\pi$.
\end{proof}

\subsubsection{Remark} The function $\la = \tau ^{-1}$ is the
spectral function of the symmetric operator $B$, as defined and studied in
Section 3 of \cite{Martin-symsamp}. Since the characteristic function $\om _B$ of $B$ is only defined
up to a unimodular constant, $\om _B$ and the phase function $\tau$ can be
chosen as described in Remark \ref{subsubsection:range} so that $\ran{\tau } = (-\infty , 0)$, $(0 , \infty )$ or $\bm{R}$.
Furthermore, if the deficiency vectors of $B$ are chosen so that $B_\alpha = B(\alpha)$, then
if $\la = \tau ^{-1}$, $\sigma (B (\alpha)) = \{  \la (\alpha + 2\pi n ) | \ n \in \bm{M} \} $
where $\bm{M} := \bm{Z} \cap \ran{\tau } = \pm \bm{N}$ or $\bm{Z}$. \label{subsubsection:taucanon}

\subsection{Subspaces with the sampling property}
\label{subsection:subsamp}

Let $\mc{H}$ be a reproducing kernel Hilbert space of functions on $\bm{R}$,
with point evaluation vectors $\delta _x$, $x\in \bm{R}$. We will call $\La := (\la _n ) _{n \in \bm{M} }$
a total orthogonal sampling sequence for $\mc{H}$ if $\la _n < \la _{n+1} $ for all $n\in \bm{M}$ and $\{ \delta _{\la _n} \} $ is a
total orthogonal set, so that any $f \in \mc{H}$ can be reconstructed from its samples $\{ f(\la _n )\}$
taken on the sampling sequence $\La$. Here $\bm{M} = \pm \bm{N}$ or $\bm{Z}$.
If $\mc{H}$ has a one-parameter family of total orthogonal sampling sequences
$\La (\alpha)$ which cover $\bm{R}$ exactly once, we will say that $\mc{H}$ has the $U(1)$ sampling
property, and if $\mc{H}$ has the sampling property and there exists an orthogonal sampling sequence $\La$ for
$\mc{H}$ which has no finite accumulation points, we will say that $\mc{H}$ has the uniformly discrete sampling property.
Reproducing kernel Hilbert spaces $\mc{H}$ of functions on $\bm{R}$ with the uniformly discrete sampling property
seem to be of greater practical value for applications such as signal processing than
those without this property. For example, suppose that one is given a RKHS $\mc{H}$ with the sampling property. One
could then attempt to use $\mc{H}$ in the same way that the space of $\Om -$bandlimited functions $B(\Om)$ is used
to sample and reconstruct continuous signals (see the Introduction). Namely, given a continuous signal $f$, \emph{e.g.}
a music signal, approximate $f$ by an element $f _\mc{H}$ of $\mc{H}$. Since a music signal is a function
of time, think of the real variable of elements of $\mc{H}$ as time. Using that $\mc{H}$ has a total orthogonal
sampling sequence $\La = (t_n ) _{n \in \bm{M}}$, one can record the samples of $f_\mc{H}$ on $\La$ to obtain the sampling
sequence $(f_\mc{H} (t_n ) )$, and then later reconstruct the approximation $f_\mc{H}$ to
$f$ from this discrete sequence. If $\mc{H}$ does not have the uniformly discrete sampling property,
the sequence of points $(t_n )$ has a finite accumulation point, and the above does not yield a practical
method for sampling and reconstructing an approximation to the continuous signal $f$.

    If $\phi$ is extreme, then as shown in \ref{subsubsection:orthog}, if
$K^2 _\phi$ has a total orthogonal sampling sequence then $\phi =F$ is an inner function.
Moreover $\sigma _e (F) \sm \{ \infty \} = \emptyset$ if and only if $F$ is a meromorphic
inner function. As observed at the beginning of Section \ref{subsection:pspec},
if $\la \notin \sigma _e (F)$ then point evaluation at $\la $ is a bounded linear functional
in $K^2 _F$. Hence if $F$ is meromorphic then $K^2 _F$ is a reproducing kernel Hilbert space of functions on $\bm{R}$.
Further suppose that  $F$ is chosen so that $M^F$ is a densely defined symmetric operator (see Theorem \ref{thm:dBdense} and use
Theorem \ref{thm:innerdB}). In this case, if the self-adjoint extensions $M ^F (\alpha )$ are defined as in Subsection \ref{subsection:pspec},
and $\sigma (M ^F (\alpha ) ) = \{ \la _n (\alpha ) \} _{n \in \bm{M}}$, ordered so that $\la _n
(\alpha ) < \la _{n+1} (\alpha )$, then it follows from the results of the previous sections
that for each $\alpha \in [0, 2\pi )$, $\sigma (M ^F (\alpha ) )$ is a strictly monotonically
increasing sequence with no finite accumulation point, the spectra of all the $M ^F (\alpha )$
cover the real line exactly once, and for each $\alpha \in [0,2\pi )$, $\{ k _{\la _n (\alpha)} ^F \}
_{n \in \bm{M}}$ is total orthogonal set of point evaluation vectors. Hence each subspace
$K^2 _F$ where $F$ is a meromorphic inner function for which $M^F$ is densely defined is a reproducing kernel Hilbert space with the
$U(1)$ uniformly discrete sampling property.

\subsubsection{Remark} It is true that spaces $K^2 _G$ for more general inner $G$ can also
have total orthogonal sets of point evaluation vectors (see Subsection \ref{subsubsection:orthog}
and \cite{Fricain}). However, if, for example, $\la \in \sigma _e (M^G) $ is an isolated point of the
essential spectrum, then it is an accumulation point of the eigenvalues
of every self-adjoint extension of $M ^G$. In particular, $\la$ will be an accumulation point of
any total orthogonal sampling sequence for $K^2 _G$. Also if $\la \in \sigma _e (G)$ then point
evaluation at $\la $ is a bounded linear functional in $K^2 _G$ if and only if the angular
derivative of $G \circ \mu ^{-1}$ exists at $\mu (\la)$. \vspace{.025truein} \\

    The following theorem applies the results of this section to provide a sufficient
condition for a subspace $\mc{H} \subset L^2 (\bm{R} , d\nu)$ to be a reproducing
kernel subspace with the $U(1)$ uniformly discrete sampling property.

\begin{thm}
    Suppose that $\nu$ is a positive Borel measure which is absolutely continuous with respect to Lebesgue measure and let $M_\nu$ be the
self-adjoint operator of multiplication by
the independent variable in $L^2 (\bm{R}, d\nu)$. Further assume that $\mc{H} \subset L^2 (\bm{R}, d\nu )$ is such
that the Cayley transform $\mu (M_\nu)$ of $M_\nu$ is a unitary dilation of its compression to $\mc{H}$  and that
$M_\nu $ has a regular simple symmetric restriction, $M ^\mc{H} _\nu$, with deficiency indices $(1,1)$ to a linear subspace
of $\mc{H}$. Then the following statements are true: \bn
\item $\mu (M_\nu )$ is the minimal unitary dilation of its compression to $\mc{H}$,
and $\nu ' (x) > 0 $ almost everywhere.
\item There is an isometric transformation $V$ which acts as multiplication by a measurable, locally $L^1$ function
which takes $\mc{H}$ onto a de Branges space of entire functions.
\item If $1/\nu '$ is a locally $L^\infty$ function and $M _\nu ^\mc{H}$ is densely defined, then
$\mc{H}$ itself is a reproducing kernel Hilbert space with the $U(1)$ uniformly discrete sampling property.
\en
\label{thm:RKHS}
\end{thm}

    Note that in condition (1), the assumption that $\nu $ is absolutely continuous to Lebesgue measure, and
that $\nu ' (x) >0 $ almost everywhere with respect to Lebesgue measure, is equivalent to
the assumption that $\nu $ is equivalent to Lebesgue measure. Recall that two measures are said to
be equivalent if they have the same sets of measure zero.

    This theorem is a strengthening of Theorem 14 of \cite{Martin-symsamp}. The proof of this theorem
will make use of the following lemma. Given a semigroup $\mf{S}$ of operators on a Hilbert space $\mc{H}$,
recall that a subspace $S \subset \mc{H}$ is said to be semi-invariant for $\mf{S}$ if $\mf{S} | _S$ is a
semigroup.

\begin{lemming}{ (Sarason \cite{Sarason})}
    Let $\mf{S}$ be a semi-group of operators on a Hilbert space $\mc{H}$. Then $S \subset \mc{H}$
is semi-invariant for $\mf{S}$ if and only if $S = S_1 \ominus S_2$ where $S_2 \subset S_1$, and
$S_1, S_2$ are invariant subspaces for $\mf{S}$. \label{lemming:semisar}
\end{lemming}

    It is easy to verify that if $S_2 \subset S_1$ are nested invariant subspaces for the semigroup $\mf{S}$,
then $S_1 \ominus S_2$ is semi-invariant for $\mf{S}$. Conversely, if $S$ is a semi-invariant subspace for $\mf{S}$,
then as in \cite{Sarason}, one can show that $\ov{\mf{S} S} =: S_1$ and $S_1 \ominus S =: S_2$ are invariant
subspaces of $\mf{S}$ satisfying $S = S_1 \ominus S_2$.

    Before beginning the proof of Theorem \ref{thm:RKHS}, it will be useful to first recall some
basic facts about unitary dilations of contractions, and to establish some notation.

\subsubsection{Unitary dilations of contractions and semigroups of contractions} \label{subsubsection:dilation}
Let $T$ be a contraction on $\mc{H}$. Recall that a unitary operator $U$ on $\mc{K} \supset \mc{H}$ is
called a unitary dilation of $T$ if $T^k = P_\mc{H} U^k | _{\mc{H}}$ for all $k \in \bm{N} \cup \{ 0 \}$.
The dilation $U$ is called minimal if $\mc{K}$ is the closure of the linear span of $U^k \mc{H} ;  \ k \in \bm{Z}$.
Any contraction $T$ has a minimal unitary dilation, and this minimal unitary dilation is unique up to
a certain natural unitary equivalence \cite[Theorem 4.3]{Paulsen}.

    If $1 \notin \sigma _p (T)$ where $\sigma _p (T)$ is the set of eigenvalues of $T$, then
the Hardy functional calculus can be used to define $T(t) := \exp \left( it \mu ^{-1} (T) \right)$
for each $t \geq 0$. The functional calculus further implies $\mf{S} := \{ T(t) \} _{t \geq 0}$ is a
semi-group with respect to multiplication (in fact a representation of ($[0, \infty) , +)$ since
$T(t)T(s) = T(s+t)  \ \forall s,t \geq 0$ and $T(0) =I$), that $\| T(t) \| \leq 1$ for all $t \geq 0$,
and that $t \mapsto T(t)$ is strongly continuous.  Any semigroup of operators on $\mc{H}$ with these
properties is called a strongly continuous one parameter contraction semigroup. Conversely, given
any such one parameter contraction semigroup, $\mf{S} = \{ T(t) \} _{t \geq 0}$, the limit
$T := \lim _{s\rightarrow 0 ^+} f_s (T (s)) $ where $f_s (z) := \frac{z-1+s}{z-1-s}$ always
exists in the strong operator topology. This limit, $T$, is a contraction on $\mc{H}$ such that $1 \notin
\sigma _p (T)$, and such that $T(t) := \exp \left( it \mu ^{-1} (T) \right)$. The contraction $T$
is called the co-generator of $\mf{S}$, and its inverse Cayley transform $\mu ^{-1} (T)$ is called
the generator. Analogously, a group $\mf{g} = \{ U(t) \} _{t \in \bm{R}}$ of operators on a Hilbert
space $\mc{K}$ is called a strongly continuous one parameter unitary group if each $U(t)$ is unitary,
if $t \mapsto U(t)$ is strongly continuous, and if $(\mf{g}, \cdot)$ is a representation of $(\bm{R},
+)$, \emph{i.e.} $U(t)U(s) = U(t+s)$, $U(0) = I$ and $U(-t) = U (t) ^{-1}$ for all $s,t \in \bm{R}$.
Such a strongly continuous unitary group of operators $\mf{g}$ on $\mc{K} \supset \mc{H}$ is called
a unitary dilation of a one parameter strongly continuous semigroup of contractions $\mf{S} = \{ T(t) \}
_{t \geq 0}$ on $\mc{H}$ if $P_\mc{H} U(t) | _\mc{H} = T(t) $ for all $t \geq 0$. Again, the dilation
is called minimal if the linear span of $U(t) \mc{H}; \ t \in \bm{R}$ is dense in $\mc{K}$, such
a minimal dilation is unique up to a natural unitary equivalence, and every strongly continuous one
parameter semigroup of contractions on $\mc{H}$ has a minimal unitary dilation \cite[Chapter III, Sect. 8-9]{Foias}.
We refer the reader to \cite{Foias} and \cite{Paulsen} for more on
basic dilation theory.

    By Stone's theorem, any strongly continuous one-parameter unitary group $\mf{g} = \{ U(t) \} _{t \in \bm{R}}$
of operators on $\mc{K}$ can be realized as $U(t) = e^{itA}$ for some densely defined and closed
self-adjoint operator $A$ acting in $\mc{K}$. The self-adjoint operator $A$ is said to generate
the unitary group $\mf{g}$. If $\mf{S} _\mc{K} = \{ U(t) \} _{t \geq 0}$, then this is
clearly a strongly continuous semigroup, and its co-generator is $\mu (A)$. Moreover, the unitary group
$\mf{g} $ is a (minimal) unitary dilation of a strongly continuous one parameter contraction semigroup
$\mf{S} = \{ T(t) \} _{s \geq 0}$ of operators on $\mc{H} \subset \mc{K}$ with co-generator $T$ if
and only if $\mu (A)$ is a (minimal) unitary dilation of $T$. In general, $\mf{g}$ is a unitary
dilation of its compression $\mf{S} _S := \{ P _S U(t) | _S \} _{t \geq 0}$ to a subspace $S \subset \mc{K}$
if and only if $S$ is semi-invariant for $\mf{S} _\mc{K}$. We will say $S$ is semi-invariant for
$\mf{g}$ if it is semi-invariant for the semi-group $\mf{S} _\mc{K}$.

    Let $M$, $M _\nu$, and $M_E$ denote the self-adjoint
operators of multiplication by the independent variable in $L^2(\bm{R})$, $L^2 (\bm{R} , d\nu)$ and
$L^2 (\bm{R} , |E(x)|^{-2} dx )$, respectively. We will use the notation $\mf{g} _M$, and $\mf{g} _\nu$
and $\mf{g} _E$ to denote the strongly continuous one-parameter unitary groups generated by $M$, $M_\nu$
and $M_E$.

    By the Beurling-Lax theorem the invariant subspaces of the semigroup $\mf{S} := \{ e^{itM} |_{H^2 (\bm{U})} \} _{t \geq 0}$
acting on $H^2 (\bm{U})$ all have the form $F H^2 (\bm{U})$, where $F$ is inner. It follows from Lemma \ref{lemming:semisar}
that each subspace $K^2 _F \subset H^2 (\bm{U})$ is semi-invariant for $\mf{g} _M$. Moreover, it is clear that $\mf{g} _M$
is the minimal unitary dilation of its compression $\mf{S} _M$ to $K^2 _F$. By Lemma \ref{lemming:isodB}, if $F=E^*/E$,
where $E$ is a de Branges function, then multiplication by $E$ is an isometry of $K^2 _F $ onto $\mc{H} (E)$. Clearly
this isometry intertwines $e^{itM}$ and $e^{itM_E}$ for all $t \in \bm{R}$. It follows that $\mf{g} _E $ is the minimal
unitary dilation of its compression, $\mf{S} _E := \{ P _{\mc{H} (E)} e^{itM_E} | _{\mc{H} (E)} \} _{t \geq 0}$, to $\mc{H} (E)$,  and that
$\mf{S} _E$ is a semi-group. Here, $P_{\mc{H} (E)}$ is the
projector of $L^2 (\bm{R} , |E(x)| ^{-2} dx)$ onto $\mc{H} (E)$.

We now proceed with the proof of Theorem \ref{thm:RKHS}.

\begin{proof}{ (of Theorem \ref{thm:RKHS}) }
    The linear map $V$ of multiplication by $\sqrt{\nu ' (x)} $ is an isometry of $L^2 (\bm{R} , d\nu)$
onto $L^2 (\bm{R})$ which maps $\mc{H}$ onto
a subspace $\mc{J} \subset L^2 (\bm{R} )$. If $M$ denotes the self-adjoint
operator of multiplication by the independent variable in $L^2 (\bm{R} )$, then $V e^{itM _\nu}
= e^{it M} V$ so that (by Remark \ref{subsubsection:dilation} above) $\mc{J}$ is semi-invariant for the semigroup $\mf{S} = \{ e^{itM} \} _{t \geq 0}$ acting
on $L^2 (\bm{R} )$.  By the Beurling-Lax theorem on invariant subspaces of this semi-group of operators
on $L^2 (\bm{R})$, and Lemma \ref{lemming:semisar},
it follows that $\mc{J} = S_1 \ominus S_2 $ where $S_2 \subset S_1$ and each $S_i$ is either equal
to $L^2 (\Om )$ where $\Om$ is a Borel subset of $\bm{R}$ or $FH^2 (\bm{U})$, where $F$ is a unimodular function.
Let $M ^\mc{J}$ denote the image of $M ^\mc{H} _\nu$ under $V$.

The subspaces $S_1$ and $S_2$ cannot both have the form $L^2 (\Om)$, as then it would follow
that $\mc{J}$ is itself invariant for $M$ so that $\mc{J} = L^2 (\La)$ where $\La$ is
a Borel subset of non-zero Lebesgue measure. It would
follow that $M ' = M | _\mc{J}$ is a self-adjoint extension of $M ^\mc{J}$, and $\sigma(M') = \La$, contradicting
the fact that the spectrum of any self-adjoint extension of $M ^\mc{J}$ is purely discrete.
Furthermore, it cannot happen that only one of either $S_1 $ or $S_2$ has the form $L^2 (\Om)$. First, if $\Om$
is a proper non-trivial Borel subset of non-zero Lebesgue measure, then $L^2 (\Om)$  neither contains nor is contained in $FH^2 (\bm{U})$ for any unimodular function $F$. Hence it cannot be that one of $S_1$, $S_2$ is equal to $L^2 (\Om)$ for such a set $\Om$
while the other is equal to $FH^2 (\bm{U})$. Suppose that either $S_1 = L^2 (\bm{R} )$ and  $S_2 = FH^2 (\bm{U})$ or $S_1 =FH^2 (\bm{U})$ and $S_2 = \{ 0 \}$. Then
multiplication by $1/F$ is an isometry from $\mc{J} = S_1 \ominus S_2$ onto $S= H^2 (\bm{L})$ or $=H^2 (\bm{U})$ respectively, and the
image $M ^S$ of $M ^\mc{J}$ under this isometry acts as multiplication by the independent variable. Consider the case
where $S= H^2 (\bm{U})$. The operator $\mu (M)| _{H^2 (\bm{U})}$ is an isometry on $H^2 (\bm{U})$. Since $1$ is not an eigenvalue of $\mu (M) | _{H^2 (\bm{U})}$, $M' = \mu ^{-1} (\mu (M)| _{H^2 (\bm{U})})$ is a densely defined symmetric operator in $H^2 (\bm{U})$. Furthermore, since
$\mu (M') = \mu (M)| _{H^2 (\bm{U})}$ is unitarily equivalent to the forward shift on $H^2 (\bm{D})$, it is
easy to see that $\dim{\dom{\mu (M') } ^\perp} = 0$ and $\dim{\ran{\mu (M')} ^\perp } =1$ so that
$\mu (M')$ and $M'$ have deficiency indices $(0,1)$. Since $\mu (M^S) = \mu (M) | _{\ran{M ^S +i} \subset H^2 (\bm{U})}$ it
follows that $\mu (M')$ is a closed isometric extension of $\mu (M^S)$. This yields a contradiction, since
$\mu (M^S)$ has deficiency indices $(1,1)$ so that the only isometric extensions of $\mu (M^S)$ have the form $\mu (M^S ) \oplus W$ on $S = \dom{\mu (M^S) } \oplus \dom{\mu (M^S)} ^\perp$ where $W$ is a rank one isometry from the one-dimensional subspace $\dom{\mu (M^S)} ^\perp$ onto the one dimensional
subspace $\ran{\mu (M ^S)} ^\perp$. All such extensions have deficiency indices $(0,0)$ and are in fact unitary. A similar
argument shows that $S$ cannot equal $H^2 (\bm{L})$ either.

    In conclusion $\mc{J} = GH^2 (\bm{U})\ominus FH^2(\bm{U})= G K^2 _{F/G}$ where $F,G$ are unimodular functions and $FH^2(\bm{U}) \subset GH^2(\bm{U})$.
Since $FH^2(\bm{U}) \subset GH^2(\bm{U})$, it follows that for any $h \in H^2(\bm{U})$ there is a $h_2 \in H^2 (\bm{U}) $ such
that $Fh = Gh _2$. That is, given any $h \in H^2(\bm{U}) $, $F/G h  \in H^2$. This implies that $F/G$ is an inner
function. Hence $\mc{J} = G \left( H^2(\bm{U}) \ominus F/G H^2 (\bm{U}) \right)$. Since elements of $\mc{J}$ have support on all of $\bm{R}$
it follows that $\nu' (x) >0$ almost everywhere $x \in \bm{R}$. Furthermore, since $\mf{g} _M$
is the minimal unitary dilation of its compression to $K^2 _{F/G}$ and hence also of its compression to $\mc{J}$, it follows that $\mf{g} _\nu$ is the minimal
unitary dilation of its compression to $\mc{H}$.

    Multiplication by $1/G$ is an isometry of $\mc{J}$ onto the subspace
$\mc{K} := K^2 _{F/G} = H^2 (\bm{U})\ominus F/G H^2 (\bm{U})$, and the image $M ^\mc{K}$ of $M ^\mc{J}$ under this isometry
is again multiplication by the independent variable. Note by Remark \ref{subsubsection:unique},
that $M ^\mc{K} = M^{K^2 _{F/G}} = M ^{F/G}$, where $M ^{F/G}$ is the symmetric linear transformation
of multiplication by $z$ in $K^2 _{F/G}$ as defined in previous sections. The transformation $M^{F/G}$ is simple
with deficiency indices $(1,1)$.
Since $M$ and hence $M ^\mc{K} = M ^{F/G}$ is regular, it follows as in the proof
of Theorem \ref{thm:dBrep} that $F/G$ is a meromorphic inner function, that there is a de Branges function $E$ such that
$F/G = E^*/E$, and that multiplication by $E$ is an isometry of $K^2 _{F/G}$ onto
$\mc{H} (E)$. In summary, if $\wt{V}$ denotes the operator of multiplication by
$v:= \sqrt{\nu ' } \frac{E}{G}$, then $\wt{V}$ is an isometry of $\mc{H}$ onto $\mc{H} (E)$
that takes $M ^\mc{H}$ onto the symmetric operator of multiplication by $z$ in $\mc{H} (E)$.

    To prove the third and final statement, let $\wt{v} $ be a member of the
equivalence class of $v$ which is bounded below on any finite interval.
Then if $k_x$ is the point evaluation vector at $x \in \bm{R}$ for $\mc{H} (E)$,
it is easy to see that $ \delta _x := (\ov{\wt{v} (x)}) ^{-1} \wt{V}^* k_x$ is such that
for any $f \in \mc{H}$, $\ip{f}{\delta _x} = f(x)$ almost everywhere. Identifying
each $f$ in $\mc{H}$ with that member of its equivalence class for which this is
true everywhere, we see that $\mc{H}$ is a reproducing kernel Hilbert space. The
fact that $\mc{H}$ has the uniformly discrete $U(1)$ sampling property follows
from the fact that any $\mc{H} (E)$ (or equivalently $K^2 _F $ where $F =E^*/E$) in which
 multiplication by $z$ is densely defined has this property.
\end{proof}

\subsubsection{Remark} In the above theorem, the assumption that $M _\nu$ has a regular simple symmetric restriction
 with deficiency indices $(1,1)$ to a linear subspace of $\mc{H}$ is equivalent to the
the assumption that the essential spectrum of the compression of $\mu (M _\nu)$ to $\mc{H}$
consists of the singleton $\{ 1 \} $.
If $\sigma _e (P _\mc{H} \mu (M _\nu) | _\mc{H})$ contains only the point $1$, then $\mc{H}$ cannot be invariant
for $M_\nu$, and hence is not reducing for $\mf{g} _\nu =\{ e^{itM _\nu} \} _{t \in \bm{R}}$. It cannot be
invariant or co-invariant for $\mf{S} := \{ e^{itM _\nu} \} _{t \geq 0}$ either, as this would imply
that it is invariant or co-invariant for $\mu (M _\nu)$. It would follow that the compression of $\mu (M _\nu)$
to $\mc{H}$ is either an isometry or a co-isometry, and hence its spectrum would be the closed unit disc. In particular,
the entire unit circle would belong to the essential spectrum of the compression of $\mu (M _\nu )$ to $\mc{H}$,
contradicting our initial assumptions. As in the proof of the theorem above, this and the
assumption that $\mc{H}$ is semi-invariant for $\mf{g} _\nu$ then implies that there is a unitary $U$ from
$\mc{H}$ onto a model subspace $K^2 _F$, where $F$ is an inner function, such that $U$ acts as multiplication
by a locally $L^1$ function. As has been demonstrated earlier, $M$ has a unique simple
symmetric restriction $M ^F$ with deficiency indices $(1,1)$ to a linear subspace of $K^2 _F$, and hence $M_\nu$ has such a restriction
$M ^\mc{H} _\nu = U^* M ^F U$ to a linear subspace of $\mc{H}$. The assumption that $\sigma _e ( P _\mc{H} \mu (M _\nu) | _\mc{H} )$
contains only the point $1$ further implies that $M ^F$ and hence $M ^\mc{H} _\nu$ must be regular.

\begin{thm}
    Let $\mc{H}$ be a reproducing kernel Hilbert space of functions on $\bm{R}$
whose reproducing kernel is positive almost everywhere with respect to Lebesgue measure,
 $\| k_x \| ^2 > 0$ a.e. $x\in \bm{R}$. Suppose that the operator of multiplication by the independent variable,
$M ^\mc{H}$ in $\mc{H}$
is a densely defined regular simple symmetric linear operator with deficiency indices $(1,1)$.
Then there is an isometry $V$ which acts as multiplication by a function which is non-zero
almost everywhere, which takes $\mc{H}$ onto a de Branges space $\mc{H} (E)$ of entire functions,
and which takes $M ^\mc{H}$ onto the symmetric operator of multiplication by $z$ in $\mc{H} (E)$.
\label{thm:RKHS2}
\end{thm}

    The proof of this theorem relies on the theory of spectral representations of
symmetric operators as developed by M.G. Krein (see for example \cite{Krein}, \cite{Silva}
or \cite[Section 2.2]{Martin-symsamp}), and tools developed in \cite{Martin-symsamp}.
For the convenience of the reader, we provide a brief summary of the background theory
needed in the proof of the above theorem.

    Let $B$ be a closed regular simple symmetric linear operator with deficiency
indices $(1,1)$ defined on a dense domain $\dom{B} \subset \mc{H}$. Let $A$ be
an arbitrary self-adjoint extension of $B$, and define the meromorphic vector-valued
function $\psi _z := (A -i ) (A- z) ^{-1} \psi _i $ where $\psi _i $
is a fixed non-zero vector in $\ker{B^* -i}$. Next define $\delta _z := \frac{\psi _z}{\ip{\psi _z}{\psi _i}}$.
By Lemma 2 of \cite{Martin-symsamp}, $\ip{\psi _x}{\psi _i } \neq 0$ for any $x \in \bm{R}$.
Furthermore, $\delta _z$ is a meromorphic vector valued
function with poles that lie off the real axis such that $\delta _x \in \ker{B^* -x}$ for each
$x \in \bm{R}$ (see for example Section 2.2 of \cite{Martin-symsamp}).
One can then define a linear map $\Phi$ of $\mc{H}$ onto a certain vector space of meromorphic functions
by $\Phi [f] (z)  := \ip{f}{\delta _{\ov{z}}}$ for any $f \in \mc{H}$. It is easy to check that the
image of $B$ under $\Phi$ acts as multiplication by the independent variable. Indeed, if $f \in \dom{B}$,
then $\Phi [Bf] (z) = \ip{Bf}{\delta _{\ov{z}}} = \ip{f}{B^* \delta _{\ov{z}}} = z \ip{f}{\delta _{\ov{z}}} =
z \Phi [f] (z)$.

    By Corollary \ref{cor:specfunk}, and Remark \ref{subsubsection:taucanon},
the deficiency vectors $\psi _\pm$ of $B$ can be chosen so that
$\sigma (B (\alpha ) )= ( \la ( \alpha + 2\pi n ) ) _{n \in \bm{M} }$ where
$\la = \tau ^{-1}$ is the spectral function of $B$ is defined on $
\ran{\tau} = (-\infty , 0 )$, $(0 , \infty )$ or $\bm{R}$ and $\bm{M} = \pm \bm{N}$ or $\bm{Z}$.
Without loss of generality, assume that $\ran{\tau} =\bm{R}$ so that $\bm{M} = \bm{Z}$.
As in Section 4 of \cite{Martin-symsamp}, one can endow the range of $\Phi $ with an inner product as follows.

 Let $\rho$ be an arbitrary positive Borel probability measure on $[0,2\pi]$, \emph{i.e} $\rho ( [0,2\pi] )
= 1$. Given any $\phi \in \mc{H}$,
\ba \ip{\phi}{\phi} & = & \int _0 ^{2\pi} \ip{\phi}{\phi} d\rho (\alpha) \nonumber \\ & = &
\int _0 ^{2\pi} \sum _{n \in \bm{Z}} \ip{\phi }{\delta _{\la  (\alpha +2\pi n )}}
\ip{\delta _{\la (\alpha +2\pi n )}}{\phi} \frac{1}{\| \delta _{\la (\alpha +2\pi n)} \| ^2} d\rho (\alpha ) \nonumber \\
& = & \sum _{n \in \bm{Z}} \int _0 ^{2\pi}   \ip{\phi }{\delta _{\la (\alpha + 2\pi n)}}
\ip{\delta _{\la (\alpha +2\pi n  )}}{\phi} \frac{1}{\| \delta _{\la (\alpha +2\pi n)} \| ^2} d \rho (\alpha)
\label{eq:Fubchange} \\ & = & \intfty \ip{\phi }{\delta _{\la(x)}} \ip{\delta _{\la(x)} }{\phi}
\frac{1}{\| \delta _{\la (x)} \| ^2 } d \rho (x) \nonumber \\  &=& \intfty
\ip{\phi }{\delta _y} \ip{\delta _y}{\phi} \frac{1}{\| \delta _y \| ^2} d\rho (\tau (y)) \nonumber \\
& =& \intfty | \hat{\phi} (y) | ^2 \frac{1}{\| \delta _y \| ^2} d\rho (\tau (y)) .\ea  In the above,
the measure $\rho$ is extended periodically to define a measure on $\bm{R}$. It follows that
$\Phi $ can be viewed as an isometry of $\mc{H}$ onto a subspace of $ L^2 (\bm{R} , d\sigma )$
where $d\sigma (x) = d\rho (\tau (x) )$. \vspace{.05truein} \\

    Choose an arbitrary $z =r e^{i\beta } \in \bm{D}$, and define
$ d\rho _z (x) = P_r (\beta -x ) dx$ where $P_r (\theta) := \frac{1-r^2}{1-2r\cos \theta +r^2}$ is the
Poisson kernel. Let $\sigma _z ' (x) dx := d \sigma _z (x) = \frac{1}{\| \delta _x \| ^{-2}} d\rho _z (\tau (x) )
= \| \delta _x \| ^{-2} P_r (\beta - \tau (x) ) \tau ' (x) dx$. Let $M _z$ denote multiplication
by the independent variable in $L^2 (\bm{R} ,d\sigma _z)$, and let $\Phi _z$ denote the map $\Phi$
viewed as an isometry from $\mc{H}$ onto $\mc{H} _z := \Phi [\mc{H} ] \subset L^2 (\bm{R} , d\sigma _z)$.
The following is taken from \cite[Theorem 11]{Martin-symsamp}.

\begin{thm}
    The subspace $\mc{H} _z \subset L^2 (\bm{R} , d\sigma _z )$ is semi-invariant for
$\mf{g} _z := \{ e^{it M_z} \} _{t \in \bm{R}}$, and $\mf{g} _z$ is the minimal unitary
dilation of its compression to $\mc{H} _z$.  \label{thm:tech}
\end{thm}

With the aid of the above facts, we are now ready to prove Theorem \ref{thm:RKHS2}.

\begin{proof}
    For each $x \in \bm{R}$ such that $k_x \neq 0$, it is clear that $k_x$ is an eigenvector
of $(M ^\mc{H}) ^* $ to eigenvalue $x$, where $M ^\mc{H}$ denotes the symmetric operator
of multiplication by the independent variable in $\mc{H}$. By the discussion preceding the proof, there is an isometry $\Phi$
of $\mc{H}$ onto a subspace $\mc{J} \subset L^2 (\bm{R} , d\sigma )$ such that
$d\sigma = \sigma ' (x) dx $, $\sigma ' , 1/\sigma '$ are locally $L^\infty $ functions,
and such that $\mc{J}$ is semi-invariant for the semi-group $\mf{S} := \{ e^{itM _\sigma} \} _{t\geq 0}$ where $M _\sigma$ denotes
multiplication by $x$ in $L^2 (\bm{R} , d\sigma )$. Furthermore, the isometry $\Phi$ can be
chosen such that if $f \in \mc{H}$, $\Phi [f] (x) = \frac{\ip{f}{\delta _x} }{\ip{\delta _i}{\delta _x}}$
for $x \in \bm{R}$ where $\delta _x \in \ker{(M ^\mc{H}) ^* -x}$. Since $\ker{(M ^\mc{H}) ^* -x} = \bm{C} \{ k_x \} $ for
each $x \in \bm{R} $ such that $k_x \neq 0$, it follows that $\delta _x = c(x) k_x$ almost everywhere
$x \in \bm{R}$ so that $\Phi [f] (x) = \frac{f(x)}{\delta _i (x) }$ almost everywhere.

In summary, $1/\sigma '$ is locally $L^\infty$, $\mc{J}$ is semi-invariant for $\mf{S} = \{ e^{it M _\sigma } \} _{t \geq 0}$,
and $M ^{\mc{J}} := M _\sigma | _{\Phi \dom{M _\mc{H} }}$ is regular simple symmetric and densely defined
with indices $(1,1)$. By Remark \ref{subsubsection:dilation} and Theorem \ref{thm:RKHS}, there is an isometry
$V$ which acts as multiplication by a
measurable function $v$ and which takes $\mc{J}$ onto a de Branges space $\mc{H} (E)$. In conclusion, multiplication
of elements of $\mc{H}$ by $\frac{ v(x) }{\phi _i (x)} $ is an isometry of $\mc{H}$
onto a de Branges space of entire functions, and which maps $M ^\mc{H}$ onto multiplication by $z$ in $\mc{H} (E)$.
\end{proof}

\begin{cor}
    In addition to the assumptions of the above theorem, suppose that $\mc{H} \subset L^2 (\bm{R} , d\nu)$, where
$\nu $ is a non-decreasing function of $x \in \bm{R}$. Let $E$ be any de Branges function such that $M _\mc{H}$, multiplication by $x$ in
$\mc{H}$, is unitarily equivalent to $M _E$, multiplication by $z$ in $\mc{H} (E)$. Let $V$ be an isometry which
takes $M ^{F_E}$ onto $M _\mc{H}$, $F_E = E^*/E$, and acts as multiplication by the function $v(x)$.
Then there exists $\phi \in B_1 (H^\infty (\bm{U}))$ such that $ \re{ \frac{1 +F _E \phi}{1-F_E \phi}}$ is the Poisson integral
of the measure $\mu$ where $d\mu (t) = | v(t) | ^2 d\nu (t)$.
\label{cor:butchs}
\end{cor}

    This corollary uses the following fact \cite[Problem 90, pg. 90]{deBranges}:

\begin{lemming}
    Let $\mc{H} (E)$ be a de Branges space, and let $\nu : \bm{R} \rightarrow \bm{R}$ be a non-decreasing function
such that \[ \intfty | f(x) / E(x) | ^2 dx = \intfty | f(x) / E(x) ^2| d\nu (x)   \] for all $f \in \mc{H} (E)$.
Then there exists a $\phi \in B_1 (H^\infty (\bm{U}))$ such that
\[ \re{ \frac{1 + F_E (z) \phi (z) }{1 -F_E (z)  \phi (z)}} = \frac{y}{\pi} \intfty \frac{d \nu (t)}{(x-t)^2 +y^2} \]
where $z = x+iy \in \bm{U}$.  \label{lemming:brangy}
\end{lemming}

\begin{proof}{ (Corollary \ref{cor:butchs})}
    Since $V$ is an isometry it follows that for any $f \in \mc{H} (E)$,
\be \| f \| _{\mc{H} (E)} ^2 = \| f/E \| ^2 _{L^2 (\bm{R} )} = \intfty \frac{ |f (t) | ^2 }{|E (t) |^2 }
| v(t) |^2 d\nu (t). \ee The corollary now follows from Lemma \ref{lemming:brangy}.
\end{proof}

\end{document}